\documentclass[a4paper,reqno]{amsart}

\usepackage{graphicx}
\usepackage{amssymb}
\usepackage{stmaryrd}
\usepackage[leqno]{amsmath}
\usepackage{enumitem}
\usepackage{stmaryrd}

\usepackage{mathrsfs}
\usepackage{mathtools}
\usepackage{exscale}
\usepackage{relsize}
\setlist[enumerate]{label=\textnormal{(\roman*)}}
\usepackage{hyperref}
\usepackage{xcolor}
\usepackage{slashed}
\usepackage{bbm}
\usepackage{bm}
\usepackage{etoolbox}

\usepackage{appendix}
\usepackage{verbatim}

\numberwithin{equation}{section}

\newcommand{\m}[1]{
\ifdefequal{#1}{1}
{\mathbbm{#1}}
{\mathbb{#1}}
}
\newcommand{\gh}[1]{\mathfrak{#1}}
\newcommand{\q}[1]{\mathcal{#1}}
\newcommand{\mc}[1]{\mathscr{#1}}
\newcommand{\bs}[1]{\boldsymbol{\mathrm{#1}}}

\newcommand\ds{\displaystyle}

\newcommand\RR{\mathbb{R}}

\newcommand\ud{\, \mathrm{d}}
\newcommand\e{\varepsilon}

\newcommand\md{\textnormal{Mod}}
\newcommand\vve{\vec \e\,}

\newcommand\vvep{\vec \e_\perp}

\newcommand\eun{\mathbf{e}_1}
\newcommand\ed{\mathbf{e}_d}

\newcommand\TD{T_\delta}
\newcommand\TS{T_*}

\newcommand\ENE{{H^1\times L^2}}

\newcommand\LL{\mathcal L}
\newcommand\BB{\mathcal G}
\DeclareMathOperator{\spn}{Span}
\DeclareMathOperator{\Id}{Id}

\renewcommand{\Im}{\text{Im}}
\renewcommand{\le}{\leqslant}
\renewcommand{\ge}{\geqslant}

\newtheorem{theorem}{Theorem}[section]
\newtheorem{corollary}[theorem]{Corollary}
\newtheorem{lemma}[theorem]{Lemma}
\newtheorem{proposition}[theorem]{Proposition}
\newtheorem{claim}[theorem]{Claim}
\theoremstyle{definition}
\newtheorem{definition}[theorem]{Definition}
\newtheorem{remark}[theorem]{Remark}

\newcommand{\haiming}[1]{{\color{purple}#1}}

\begin{document}

\title[Construction of symmetric multi-solitary waves]{Construction of multi solitary waves with symmetry for the nonlinear damped Klein-Gordon equation }

\author[Raphaël Côte]{Raphaël Côte}
\address{IRMA UMR 7501, Université de Strasbourg, CNRS, F-67000 Strasbourg, France \newline
\& USIAS, Université de Strasbourg, F-67083 Strasbourg France}
\email{rcote@unistra.fr}

\author[Haiming Du]{Haiming Du}
\address{Beijing Institute of Mathematical Sciences and Applications, Beijing 101408, China}
\email{duhaiming@bimsa.cn}


\thanks{RC acknowledges support from the University of Strasbourg Institute for Advanced Study (USIAS) for a Fellowship within the French national program ``Investment for the future'' (IdEx-Unistra), and from the ANR Project MOSICOF ANR-21-CE40-0004.}
		
\begin{abstract}
We are interested in the nonlinear damped Klein-Gordon equation 
\[ \partial_t^2 u+2\alpha \partial_t u-\Delta u+u-|u|^{p-1}u=0 \]
on $\RR^d$ for $2\le d\le 5$ and energy sub-critical exponents $2 < p < \frac{d+2}{d-2}$.
   
We construct multi-soliton, that is, solutions which behave for large times as a sum of decoupled solitons, in various configurations with symmetry:  this includes multi-soliton whose soliton centers lie at the vertices of an expanding regular polygon (with or without a center), of a regular polyhedron (with a center), of a higher dimensional regular polytope, or on a line. We give a precise description of these multi-solitons, in particular the interaction between nearest neighbour solitons is asymptotic to $\ln (t) - \frac{d-1}{2} \ln \ln t$ as $t \to +\infty$.

We also prove that in any multi-soliton, the solitons cannot share all the same sign. 

Both statements generalize and refine the results of \cite{F98}, \cite{Nak} and are based on the analysis developed in \cite{CMYZ,CMY}.
\end{abstract}

\maketitle

\section{Introduction}
 
\subsection{Setting of the problem}

We consider the nonlinear focusing damped Klein-Gordon equation
\begin{equation}\label{DLKG}\tag{DKG}
    \left\{\begin{array}{l}
    \partial_t^2 u+2\alpha \partial_t u-\Delta u+u - f(u) =0, \\
    \left(u,  \partial_t u\right)|_{t=0}=\left(u_0(x), u_1(x)\right),
    \end{array} \quad(t, x) \in \mathbb{R} \times \mathbb{R}^d\right.
\end{equation}
where  $\alpha>0$, $1\le d\le 5$, and $f(u) = |u|^{p-1}u$. We will always assume that the exponent $p$ lies in the energy sub-critical range: $2< p <\infty$ for $d=1,2$ and $2<p<\frac{d+2}{d-2}$ for $d=3,4,5$. There is a natural  energy functional related to this equation:
\begin{equation}\label{energy}
    E(\vec u)=\frac 12 \int_{\RR^d} \big( |\nabla u|^2 + u^2 + (\partial_t u)^2 - 2 F(u) \big) \ud x, \end{equation}
where we denoted
 \[ F(u) = \int_0^u f(x) \ud x = \frac{1}{p+1} |u|^{p+1}. \]
One can rewrite \eqref{DLKG} as a first order system for $\vec u=(u, \partial_t u)=(u,v)$
\begin{equation*}
    \left\{\begin{aligned}
    \partial_t u & = v\\
    \partial_{t} v & = \Delta u - u + f(u) - 2 \alpha v.
    \end{aligned}\right.
\end{equation*}
It follows from \cite{BRS} that the Cauchy problem for~\eqref{DLKG} is locally well-posed in the energy space
$ H^1\left(\mathbb{R}^d\right) \times L^2\left(\mathbb{R}^d\right). $

\begin{proposition}[Local well-posedness, {\cite[Theorem 2.3]{BRS}}] \label{prop:lwp}
Let $d\le 6$.
For every data $\vec u_0 = (u_0,u_1)$ in $H^1(\m R^d) \times L^2(\m R^d)$, there exists a unique maximal mild solution of \eqref{DLKG} 
\[
		\vec u = (u,\partial_t u) \in \mc C\left([0, T_+), H^1\left(\mathbb{R}^d\right) \times L^2\left(\mathbb{R}^d\right)\right),
\]
		where the maximal time of existence $T_+$ is bounded below by a function of the norm of the initial data $\left\|u(0),\partial_{t}u(0)\right\|_{H^1\times L^2}$ only. 
		Moreover, if $3 \leqslant d \leqslant 6$,
$u \in L_t^{\frac{d+2}{d-2}}L_x^{\frac{2(d+2)}{d-2}}([0,T]\times \RR^d)$ for all $T <T_+$
	and furthermore, the following properties hold.
\begin{itemize}
    
    \item For $T <T_+$, the flow map
    \[  \begin{array}{r@{\ }c@{\ }l}  H^1(\m R^d) \times L^2(\m R^d) & \to & \mc C([0,T], H^1(\m R^d) \times L^2(\m R^d)) \\
     \vec v_0 & \mapsto & \vec{v} = (v, \partial v) \text{ solution to } \eqref{DLKG} \text{ with data } \vec v_0
     \end{array} \]
    is well defined and $\mc C^1$ on a neighbourhood of $\vec u_0$.


    \item If $T_+<\infty$, then 
    \[ \limsup _{t \rightarrow T_+}\|\vec{u}(t)\|_{H^1(\m R^d) \times L^2(\m R^d)}=+\infty. \]
    \item The energy of a solution decreases with time, and more precisely,
    \begin{equation}
		\forall t_1, t_2 \in [0,T_+), \quad E(\vec u(t_2))-E(\vec u(t_1)) = -2 \alpha \int_{t_1}^{t_2} \|\partial_t u(t)\|_{L^2}^2 \ud t.
    \end{equation}
\end{itemize}
	\end{proposition}
	
	In this paper, we are interested in constructing multi-solitary wave solutions, and in giving detailed description of these solutions. These are related to the ground state $Q$, which is the unique positive and radial $H^1(\m R^d)$ solution of
	\begin{equation}\label{eq:elliptic}
	-\Delta Q + Q - f(Q)=0, \quad x\in \RR^d.
	\end{equation}
	(See~\cite{BL,K}.)
	The ground state generates the stationary solution of~\eqref{DLKG},
 \[ \vec Q=(Q,0). \]
	The function $-\vec Q$ as well as any translate $\vec Q(\cdot-z_0)$ are also solutions of~\eqref{DLKG}, which we call solitons.
	
	\bigskip

	A multi-solitary wave is a solution to \eqref{DLKG} defined  for large positive times, and behaving as a sum of decoupled ground states as $t \to +\infty$ (that is the relative distance of the centers of the solitons tends to $+\infty$). The question of the existence of multi-solitary waves for~\eqref{DLKG} was first addressed by Feireisl in~\cite[Theorem 1.1]{F98}: under suitable conditions on $d$ and $p$, for an even number of solitary waves on a regular (planar) polygon, with alternating signum of nearest neighbours.
	His construction is based on variational and symmetry arguments to treat the instability direction of the solitary waves.

\bigskip

 Our objective is to present additional examples of multi-solitary waves  and describe precisely the asymptotics, in the context of configurations with symmetry. We also prove that multi-solitons necessarily have solitons of opposite signs: this generalizes a statement of \cite{CMYZ} for 2-solitons to the case of any numbers of solitons.

\subsection{Previous Results}

The investigation of the long-term asymptotic behavior of solutions to the damped Klein-Gordon equation in relation to bound states has been tackled in several articles; see \emph{e.g.}~\cite{BRS,F94,F98,J,LZ}. Significantly, under some conditions on $d$ and $p$, results in~\cite{F98,LZ} state that for any sequence of time, any global bounded solution of~\eqref{DLKG} converges to a sum of decoupled bound states after extraction of a subsequence of times.

In~\cite{BRS}, for radial solutions in dimension $d\geq 2$, the convergence of any global solution to one equilibrium is proved to hold for the whole sequence of time. As discussed in~\cite{BRS}, such results are closely related to the general soliton resolution conjecture for global bounded solutions of dispersive problems; see~\cite{DKM,DJKM} for details and results related to this conjecture for the undamped energy critical wave equation. 

\bigskip

The existence and properties of multi-solitary waves is a classical question for integrable models (see for instance~\cite{Miura} for the Korteweg-de Vries equation and \cite{ZS} for the $1$D cubic Schr\"odinger equation). Since then the same question has been addressed for various non-integrable and undamped nonlinear dispersive equations: we refer for example  to \cite{Mar05,MM06,CL11,CMMgn,CMkg} (and the references therein)  for the generalized Korteweg-de Vries equation, the nonlinear Schr\"odinger equation, the Klein-Gordon equation and the wave equation, in situations where stable or unstable ground states are involved. In those works, the distance between two travelling waves is asymptotic to $Ct$ for some $C>0$, as $t\to\infty$, and so interactions decay exponentially fast in time. There also exist multi-solitary wave with strong interaction (decaying only polynomially), where the distance between solitons is logarithmic: see for example~\cite{Jkdv,MN,TVNkdv,TVN} for the Schr\"odinger and Korteweg-de Vries  type equations and systems. Regarding strong interactions, the early reference is actually the seminal paper \cite{MMwave}: it is concerned with the energy critical wave equation, for which the distance is in $Ct$ but solitons decay only algebraically, see also \cite{MR18} on the nonlinear Schrödinger equation.

Note that the logarithmic distance in the latter works is exceptional while it is the generic behavior for the damped equation \eqref{DLKG}.

See also~\cite{JJnon,JJ,JL} for works on the non-existence, existence and classification of radial two-bubble solutions for the energy critical wave equation in large dimensions. The construction of (center-)stable manifolds in the neighbourhood of unstable ground state was addressed in several situations, see \emph{e.g.} \cite{BJ,KS,KNS,MaMeNaRa,NS1}. While the initial motivation and several technical tools originate from some of the above mentioned papers, we point out that the present article is essentially self-contained, except for the local Cauchy theory for~\eqref{DLKG} (see \cite{BRS}) and elliptic theory for \eqref{eq:elliptic} and its linearization (we refer to \cite{BL,CMkg,K}).

\subsection{Main results}

As our construction will enjoy some symmetry, let us introduce the context for it.
	
Let $G$ be a subgroup of $O_d(\m R)$ and $\tau: G \to \{ \pm 1 \}$ be a homomorphism. $G$ acts on $\vec v \in H^1(\m R^d) \times L^2(\m R^d)$ via
\[ \forall R \in G, \quad R . \vec v(x) = \vec v(R^{-1} x), \]
so that $(R R') \vec v = R . (R' .\vec v)$.

We will work in the energy space with symmetry
\[ \q H_\tau = \{ \vec v \in H^1(\m R^d) \times L^2(\m R^d) : \forall R \in G, \ \forall x \in \m R^d, \quad R . \vec v = \tau(R) \vec v \} \]
Due to the uniqueness statement for the Cauchy problem in Proposition \ref{prop:lwp}, and as the equation \eqref{DLKG} is invariant under $O_d(\m R)$ and sign changes, if $(u_0,u_1) \in \q H_\tau$ then the solution $\vec u(t)$ to \eqref{DLKG} belongs to $\q H_\tau$ for all times $t$ where it is defined. Indeed, for any $R \in G$, $R. \vec u(t)$ and $\tau(R) \vec u(t)$ are both solution to \eqref{DLKG} with the same initial data $R . (u_0,u_1) = \tau(R) (u_0,u_1)$, and so coincide.

\bigskip

We denote $\Omega$ a finite set of $\m R^d$, which is meant to denote the centers of the solitons. It is also natural to ask that $\Omega$ enjoys some symmetry; actually we require a more robust property than a simple action of the group $G$ on these points.

\begin{definition} 
We say that a configuration $(\Omega, G)$ is \emph{rigid} if $G$ is a subgroup of $O_d(\m R)$ acting on $\Omega$, with transitive action over $\Omega \setminus \{0 \}$ ($\Omega$ represents the vertices of a polytope, and possibly its center), and they satisfy the following property. Let $\Omega'$ be another finite subset of $\m R^d$ with the same cardinality as $\Omega$, and on which $G$ acts as well. Then there exist $\lambda > 0$, $\tau \in \m R^d$ and $R \in O_d(\m R)$, and a (one-to-one) map $\phi: \Omega \to \Omega'$ such that
\[ \forall \omega \in \Omega, \quad \phi(\omega) = \lambda R \omega + \tau. \]
Notice that $R$ can be chosen to depend smoothly on $\Omega'$.
\end{definition}

 \bigskip
 
Together with the centers, we must specify the sign of each soliton. For this, let $\sigma: \Omega \to \{ \pm 1 \}$ be equivariant in the sense that there exists a homomorphism $\tau: G \to \{ \pm 1 \}$ such that
\begin{equation} \label{def:equiv_ts}
\forall R \in G, \ \forall \omega \in \Omega, \quad \sigma (R\omega) = \tau(R) \sigma(\omega).
\end{equation} 
This condition ensures that
\begin{equation} \label{eq:S_sym} 
\sum_{\omega \in \Omega} \sigma(\omega) (Q(\cdot - \omega),0) \in \q H_\tau
\end{equation}
has symmetry.

\bigskip
 
Our first main result is that under appropriate symmetry conditions on $\Omega$, we are able to construct a multi-soliton whose centers are placed at a dilatation of $\Omega$ (which expands with time $t$ essentially like $\ln t$), and we can give a fine description of its behavior. The precise statement is as follows.

\begin{theorem} \label{thm:abs}
Let $(\Omega, G)$ be a rigid configuration, and $\sigma: \Omega \to \{ \pm 1 \}$ equivariant in the sense of \eqref{def:equiv_ts}.

We assume that $\Omega$ has its center of mass $0$, that is $\sum_{\varpi \in \Omega} \varpi =0$, and that there exists $\gamma >0$ so that for all $\omega \in \Omega$, denoting $\Omega_\omega$ the subset of nearest neighbours of $\omega$
\begin{equation}\label{symmetry}
\sigma(\omega) \sum_{\varpi \in \Omega_\omega} \sigma(\varpi) (\omega-\varpi) = \gamma \omega.
\end{equation} 
Then, there exist $\lambda_\Omega >0$, $c_\Omega \in \m R$ and a $|\Omega|$-soliton $\vec u \in \mc C([0,+\infty), \q H_\tau)$ of the form
\begin{equation} \label{eq:multi}
u(t,x) = \sum_{\omega \in \Omega} \sigma(\omega) Q \left(x-d(t) \omega \right) + \e(t,x),
\end{equation}
such that $\| \e(t) \|_{H^1} + \| \partial_t u(t) \|_{L^2} = O(t^{-1})$ and 
\begin{equation}\label{dis:new}
d(t) = \lambda_{\Omega} \left(\ln t - \frac{d-1}{2} \ln\ln t \right) + c_{\Omega}+ O \left( \frac{\ln \ln t}{\ln t} \right).
\end{equation}
\end{theorem}

\begin{remark}
The assumption that the center of mass of $\Omega$ is $0$ can be relaxed by translation invariance; we made the choice for simplicity in the notations. 

The assumption \eqref{symmetry} is the key dynamical assumption which sustains the ansatz of the multi-soliton solution whose centers keep a shape expanding in time.

$\lambda_\Omega$ is the inverse of the minimal distance between two vertices see \eqref{def:lambda_O}, $c_\Omega$ is explicit (but less geometric), a definition is given at the end of the proof of the theorem.
\end{remark}

\begin{corollary} \label{cor:sym}
The regular polytopes together with their symmetry groups  in dimension $n \le d$ are rigid in $O_d(\m R)$.

As a consequence, there exist multi-solitons (in the sense of Theorem \ref{thm:abs}) in the following cases:
\begin{itemize}
\item $\Omega$ is a polyhedron with a center, with $d \ge |\Omega| +1$ and $\sigma(0) = -1$, $\sigma(\omega) = 1$ for $\omega \in \Omega \setminus \{ 0 \}$, $\gamma =1,\quad \lambda_\Omega^{-1} =|\omega|$, in the cases of
\begin{enumerate}
\item (2D) an equilateral triangle, a square, or a regular pentagon;
\item (3D) a regular tetrahedron ($ |\Omega| = 5$), octahedron ($ |\Omega| = 7$), cube ($ |\Omega| = 9$), or icosahedron ($ |\Omega| = 12$);
\item a regular simplex;
\item a regular orthoplex (dual of the hypercube);
\end{enumerate}
\item a $2K$ regular polygon (2D) with alternating signs (without center): viewing $\Omega$ lying on a plane which we identify with $\m C$, $\Omega = \{ e^{i k \pi/K} : k \in \llbracket 0, 2K-1 \rrbracket \}$ and $\sigma(e^{i k \pi/K}) = (-1)^k$. Then $\ds \gamma=4\sin^2\left(\frac{\pi}{|\Omega|} \right)$, $\ds \lambda_{\Omega}^{-1}= 2|\omega|\sin\left( \frac{\pi}{|\Omega|}\right) $.
\item a hypercube of any dimension and alternating signs (without the center): for some $n \le d$, $\Omega = \{ \pm 1 \}^n \times \{ 0 \}^{d-n}$,  $|\Omega| = 2^n$, $\gamma=2$, $\lambda_{\Omega}^{-1}=\frac{2}{\sqrt{n}}|\omega|$, and for $(\iota_1, \ldots, \iota_n) \in \{ \pm 1 \}^n$,
\[ \sigma (\iota_1, \ldots, \iota_n, 0, \dots, 0) = \prod_{k=1}^n \iota_k. \]
\item  a hypercube in dimension at least
$5$ with the center and alternating signs: let $5\le n \le d$, $\Omega = \{ \pm 1 \}^n \times \{ 0 \}^{d-n} \cup \{(0, \dots, 0) \}$,  $|\Omega| = 2^n+1$, $\gamma=2$, $\lambda_{\Omega}^{-1}=\frac{2}{\sqrt{n}}|\omega|$, and for $(\iota_1, \ldots, \iota_n) \in \{ \pm 1 \}^n$,
\[ \sigma (\iota_1, \ldots, \iota_n, 0, \dots, 0) = \prod_{k=1}^n \iota_k, \quad \sigma(0, \dots, 0) \in \{ \pm 1 \}. \]
\end{itemize}
\end{corollary}

Let us also mention a case where the symmetry group $G$ is not the full symmetry group of $\Omega$ (but still rigid), and which corresponds to the case of the orthoplex without a center, in even dimension.

\begin{corollary} \label{cor:orthoplex}
Let $d \ge 2$ be \emph{even}. There exist multi-solitons (in the sense of Theorem \ref{thm:abs}) where $\Omega = \{ \pm {\bf e}_j : j=1, \dots, d\}$ is the orthoplex ($(\eun,\ldots,\ed)$ is the canonical base of $\m R^d$), $\sigma({\bf e}_j) = \sigma(- {\bf e}_j) = (-1)^j$ and
\[ G = \{ g \in O_d(\m R) : \exists \epsilon(g) \in \{ \pm 1\}, \forall \omega \in \Omega, g(\omega) \in \Omega \text{ and } \sigma(g(\omega)) = \epsilon(g) \sigma(\omega) \}. \]
In that case, $\gamma = 2$ and $\lambda_\Omega = \frac{1}{\sqrt 2}$.
\end{corollary}

\begin{remark}
We emphasize that the regular (2D) polygons with at least 6 vertices and a center are not covered by the results of Corollary \ref{cor:sym}, because they do not meet the condition \eqref{symmetry} (for the hexagon, $\gamma=0$ and for higher number of vertices, the nearest neighbour of a vertex is not the center so that $\gamma <0$). We conjecture that when there are at least 7 vertices, no such multi-soliton exists; the case of the hexagon is critical.

It is interesting to note that actually, there exist \emph{stationary} (sign-changing) solutions whose mass is concentrated near the vertices of a regular polygon with 7 vertices or more: they are constructed in \cite{MPW12}.

Regarding 3D polyhedra, the case of the regular dodecahedron with a center is not covered for the same reason: the distance to the center over edge length ratio is $\frac{\sqrt 3 (1+\sqrt 5)}{4} >1$, so that the center is not the nearest neighbour and $\gamma <0$.

The octahedron (or the orthoplex in dimension at least 4) without a center is not treated either: the equivariance condition would imply that all the signs are equal (because for any pair of vertices, there is a rotation $R$ which maps one to the other while fixing two other vertices so that $\tau(R)=1$); but in this case, the expanding dynamics condition \eqref{symmetry} is not satisfied (or simply apply Theorem \ref{thm:non-existence} below).

This is related to the fact that the graph of the edges of the octahedron is not bipartite. The same reason applies to the case of the dodecahedron (with or without center) and to the icosahedron without a center, or to the other regular polytopes in dimension at least 4.

The case of a hypercube with a center (and alternating signs on the vertices) is possible in dimension $d \ge 5$ because the center is not a nearest neighbour of any vertex (distance $\sqrt d > 2$), and so the dynamics is essentially that of the hypercube without a center. The hypercube in dimension 4 is critical.

Finally, regarding the orthoplex without a center in odd dimension, a geometric setting as in Corollary \ref{cor:orthoplex} is possible, but then \eqref{symmetry} holds with $\gamma=0$, and this is not covered by Theorem \ref{thm:abs}.
\end{remark}

\begin{remark}
\cite{F98} considers the case of the $2K$ regular polygon in 2D, and \cite{Nak} considers the orthoplex with a center. Regarding these earlier works, we provide a more refined description given in \eqref{dis:new}. 
\end{remark}

\begin{remark}
One simple example of a non rigid configuration is the case of a rectangle (with $G =D_2 \simeq (\m Z/2\m Z)^2$ the dihedral group of order $4$). The rectangle with a center and aspect ratio $1:a$ with $1 < a < \sqrt 3$ (so that the center is still the nearest neighbour of the vertices) satisfies all the other conditions, but we cannot ensure that the aspect ratio remains constant over time only using symmetry argument; in other words, that the configuration is rigid (hence the name).
\end{remark}

The construction of strongly interacting multi-solitons is always delicate, as they are unstable objects. The symmetry assumptions here are done so that the dynamics is led by a single real function (the distance to the center $d$). If more freedom is allowed, one is required to understand an unstable system of ODE under perturbations, which can be untractable. The case of solitons on a line for example in $\m R^d$, $d \ge 2$, is not an immediate extension of the case $d=1$ (treated in \cite{CMYZ}), and is the object of Theorem \ref{th:2} below.

We devoted some effort to abstract the hypothesis needed on the symmetry group $G$ to be amenable to the analysis, keeping in mind we were interested mainly in configuration on regular polytopes as stated in Corollary \ref{cor:sym}. One point is that the multi-soliton may lie in a strict subspace of the  ambient space (which can therefore be higher dimensional), so that the action of $G$ is in particular not irreducible: the rigidity condition on $(\Omega, G)$ appears naturally in that context.

\bigskip

Our second construction result is focused on multi-soliton with center on a line. Here is the precise statement.

\begin{theorem} \label{th:2}
For any $K\ge 2$, $\sigma \in \{ \pm 1 \}$ and direction $e \in \m S^{d-1}$, there exists a global solution $\vec u \in \mc C([0,+\infty), H^1 \times L^2)$ of~\eqref{DLKG} such that for all $t\in [0,\infty)$,
\begin{equation}\label{eq:th:2}
\left\|u(t)-\sigma \sum_{k=1}^K(-1)^kQ(\cdot- z_k(t)e)\right\|_{H^1}
+\|\partial_{t}u(t)\|_{L^{2}}\lesssim t^{-1},
\end{equation}
and for all $k=1,\dots, K$,  
\begin{equation} \label{eq:th:z}
z_k(t) =  \left( k - \frac{K+1}{2} \right) \left( \ln t - \frac{d-1}{2} \ln \ln t \right) +  O ( 1),
\end{equation}
as $t\to \infty$.
\end{theorem}

It is an analogue of Theorem 1.3 in \cite{CMY} which was done in dimension 1, and with only a lesser control on the $z_k$ (only $O(1)$ instead of convergence, see also remarks \ref{nb:stab_ODE} and \ref{rem:last}). We emphasize that the proof requires an additional involved argument with respect to Theorem \ref{thm:abs}. Indeed, the case of the line corresponds to a cylindrical symmetry around the axis $\m R e$: all the elements of this axis are therefore fixed point, and so the configuration of the soliton centers, which lie on the axis, is \emph{not rigid}. As a consequence, one has to study the ODE \emph{system} for the solitons' center dynamics, and not only a scalar ODE. Now, this system is similar to, but not the same as the 1D one, so we have to make sure that it enjoys the suitable stability properties. All details and proofs are given in Section \ref{sec:line} and the Appendix \ref{app}.

\bigskip

Observe that in all the above constructions, the solitons never share a common sign. Our last result is that this feature is general: any multi-soliton necessarily has solitons of opposite sign. First we precise which definition of a multi-soliton we use: in the spirit of \cite{CMYZ,CMY}, in the context of \eqref{DLKG}, convergence on a sequence of times is sufficiently strong a requirement.

\begin{definition}\label{def-multi wave}
A $K$-soliton $\vec{u}$ is a solution to \eqref{DLKG} defined for sufficiently large times, i.e. $\vec u \in \mc{C}\left([T, \infty), H^1 \times L^2\right)$  for some $T \in \mathbb{R}$, and such that there exists  a sequence of times $t_n \rightarrow \infty$ and for each $k = 1, \dots, K$, a sign $\sigma_k \in \left\lbrace  \pm 1\right\rbrace $ and a sequence of points  $( \xi_{n,k})_{n \in \m N}$  in $\mathbb{R}^{d}$  such that
\[
	\left\|u\left(t_n\right)-\sum_{k=1}^K \sigma_k Q\left(\cdot-\xi_{n,k}\right)\right\|_{H^1}+\left\|\partial_t u\left(t_n\right)\right\|_{L^2} \to 0 \quad \text{as } n \to +\infty,
\]
and for each $i,j \in \llbracket 1, K \rrbracket$ with $i \ne j$, $\left|\xi_{i,n}-\xi_{j,n}\right| \rightarrow \infty$ as $n \to +\infty$.	
\end{definition}

\begin{theorem}[Non existence of same sign multi-soliton]\label{thm:non-existence}
Let $K \ge 2$, and let $\vec u$ be a $K$-soliton in the sense of Definition \ref{def-multi wave} with signum $\sigma_k \in \{\pm 1 \}$ for $k=1, \dots,K$. Then there exists $1 \le i,j \le K$  such that $\sigma_i \ne \sigma_j$.
\end{theorem}

The result is intuitively very clear: if all solitons have same sign, at leading order the interaction is attractive, and so the various solitons cannot spread away. However, the proof is more involved. Besides the fact that, in our definition, the multi-soliton configuration is only attained on a sequence of time, there is a major roadblock in the informal argument. The dynamics of the interaction of the center of mass alone does show attraction: for example, one can prove for this unperturbed system that the maximum distance between soliton decreases over time. This, however, is not the case for the minimal distance, and there are examples where the minimal distance actually increases on some interval of time. The main obstacle is that the perturbation size is actually related to the \emph{minimal} distance between solitons, and so it can be leading order in the dynamics: but then, there is no monotonicity.

Instead, we show a conditional stability result, and prove that the multi-soliton configuration holds for all times, with some quantitative bounds, via a bootstrap argument with decay in time. The decay rate however is not compatible with energy expansion, which is a contradiction.

\subsection{Organisation of the paper}

In Section \ref{sec:dyn}, we give general estimates regarding the dynamics at leading order of a sum of decoupled solitary waves. Using those, we prove Theorem \ref{thm:abs} and its Corollary \ref{cor:orthoplex} in Section \ref{sec:const}; Theorem \ref{thm:non-existence} in Section \ref{sec:non}; and Theorem \ref{th:2} in Section \ref{sec:line}. In the Appendix \ref{app}, we provide a stability analysis for the ODE system of the solitons' centers when they lie on a line, and under an amenable perturbation.

\subsection{Acknowledgment}

We would like to deeply thank the anonymous referee for very insightful comments which improved this article, and to which we are indebted for Corollary \ref{cor:orthoplex}.

\section{Dynamics of multi-solitary waves}  \label{sec:dyn}

In this section, we prove general results on solutions of \eqref{DLKG} close to the sum of decoupled solitary waves.

\subsection{Some notations and properties of ground states}

First let us introduce a few basic notations.
 Let $\{\eun,\ldots,\ed\}$ denote the canonical basis of $\RR^d$.
We denote by $\q B_d(\rho)$ (respectively, $\m S_{d-1} (\rho)$) the closed ball (respectively, the sphere) of $\m R^d$ of center the origin and of radius $\rho>0$, for
the usual norm $|\xi|= \sqrt{\sum_{j=1}^d \xi_j^2}$ (and simply $\q B_d$ and $\m S_{d-1}$ the closed unit ball and the unit sphere, when $\rho=1$).

We denote by $\q B_\ENE(\rho)$ the ball of $H^1\times L^2$ of center the origin and of radius $\rho>0$
for the norm 
\[ \|\big( \begin{smallmatrix} \e\\ \eta \end{smallmatrix}\big)\|_\ENE
=\big(\| \e\|_{H^1}^2+\| \eta\|_{L^2}^2\big)^{1/2}. \]
We denote $\langle \cdot, \cdot \rangle$ the $L^2$ scalar product for real valued functions $u_i$ or vector-valued functions $\vec u_i = (u_i,v_i)$ ($i=1,2$):
\[ \langle u_1, u_2 \rangle :=  \int u_1(x) u_2(x) \ud x, \quad \langle \vec u_1, \vec u_2 \rangle :=  \int u_1(x) u_2(x) \ud x +\int v_1(x) v_2(x) \ud x. \]

It is well-known that the operator
\begin{equation*}
\LL = -\Delta+1-p Q^{p-1}
\end{equation*}
appearing after linearization of equation~\eqref{DLKG} around $\vec Q$,
has a unique negative eigenvalue $- \nu_0^2$ ($\nu_0>0$). We denote by $Y$ the corresponding normalized eigenfunction
(see Lemma~\ref{le:L} for details and references).
In particular, it follows from explicit computations that setting
\begin{equation*}
\nu^\pm = - \alpha \pm \sqrt{\alpha^2+\nu_0^2}\quad\text{and}\quad
\vec Y^\pm = \begin{pmatrix}
Y\\ \nu^\pm Y
\end{pmatrix},
\end{equation*}
the function $\vve^{\pm}(t,x) = \exp(\nu^\pm t) \vec Y^\pm(x)$ is the solution of the linearized problem
\begin{equation}\label{eq:lin}
\left\{\begin{aligned}
\partial_t \e & = \eta\\
\partial_t \eta & = -\LL \e - 2 \alpha \eta.
\end{aligned}\right.
\end{equation}
Since $\nu^+>0$, the solution $\vve^+$ illustrates the exponential instability
of the solitary wave in positive time.
In particular, we see that the presence of the damping $\alpha>0$ in the equation
does not remove the exponential instability of the Klein-Gordon solitary wave.
An equivalent formulation of instability is obtained by setting
\begin{equation*}
\zeta^\pm = \alpha \pm \sqrt{\alpha^2+\nu_0^2}
\quad\text{and}\quad
\vec Z^\pm = \begin{pmatrix}
\zeta^\pm Y\\ Y
\end{pmatrix}
\end{equation*}
and observing that for any solution $\vve$ of~\eqref{eq:lin},
\begin{equation*}
a^\pm = \langle \vve , \vec Z^\pm\rangle \quad \text{satisfies} \quad
\frac{\ud a^\pm} {\ud t} = \nu^\pm a^\pm.
\end{equation*}

We denote the ground state $Q$ by $Q(x)=q\left(|x|\right)$ where $q>0$ satisfies
\begin{equation} \label{def:q}
q''+\frac{d-1}{r}q'-q+q^{p}=0,\quad q'(0)=0,\quad \lim_{r\to \infty}q(r)=0.
\end{equation}
It is well-known and easily checked that for a constant $\kappa>0$,
for all $r>1$,
\begin{equation}\label{Qdec}
\left|q(r)-\kappa r^{-\frac{d-1}{2}}e^{-r}\right|+\left|q'(r)+\kappa r^{-\frac{d-1}{2}}e^{-r}\right|
\lesssim r^{-\frac{d+1}{2}}e^{-r}.
\end{equation}
Due to the radial symmetry, there holds the following cancellation (which we will use repetitively):
\begin{equation} \label{eq:Q_sym}
\forall i \ne j, \quad \int \partial_{x_i} Q(x) \partial_{x_j} Q(x) \ud x =0.
\end{equation}
Let
\[
\LL = -\Delta+1-p Q^{p-1} ,\quad
\langle \LL \e,\e\rangle = \int \Big(|\nabla \e|^2+\e^{2} - p  Q^{p-1} \e^2\Big) \ud x.
\]
We recall standard properties of the operator $\LL $ (see \emph{e.g.}~\cite[Lemma 1]{CMkg}).

\begin{lemma}\label{le:L}\cite[lemma 1]{CMkg}
\emph{1. Spectral properties.} The unbounded operator $\LL $ on $L^2$ with domain $H^2$ is self-adjoint, its continuous spectrum is $[1,\infty)$, its kernel is $\spn\{\partial_{x_j}Q : j= 1,\ldots,d\}$ and it has a unique negative eigenvalue $-\nu_{0}^{2}$, with corresponding smooth normalized radial eigenfunction $Y$ $(\|Y\|_{L^2}=1)$
		Moreover, on $\RR^d$,
		\begin{equation*}
		\left|\partial^{\beta}_{x}Y(x)\right|\lesssim e^{-\sqrt{1+\nu_{0}^{2}}\left|x\right|}\quad \text{for any } \beta=\left(\beta_1 ,\ldots,\beta_d \right)\in \mathbb{N}^d.
		\end{equation*}

\emph{2. Coercivity property.} There exists $c>0$ such that, for all $\e\in H^{1}$,
		\begin{equation} \label{est:coer_L}
		\langle \LL \e,\e\rangle\ge c
		\|\e\|_{H^{1}}^{2}-c^{-1}
		\bigg(\langle \e,Y\rangle^{2} + \sum_{j=1}^d \langle \e,\partial_{x_j}Q\rangle^{2}\bigg).
		\end{equation}

\end{lemma}

\subsection{Static results}

Let $K \ge 1$. We will often use bold symbols to denote families indexed by $1, \dots, K$: let $\bm z = (z_k)_{1 \le k \le K} \in (\m R^d)^K$ be centers of mass, $\bm \ell = (\ell_{k} )_{1 \le k \le K} \in (\m R^d)^K$ be Lorentz parameters and $\bm \sigma = (\sigma_{k} )_{1 \le k \le K} \in \{ \pm 1 \}^K$ be signum (in this subsection 2.2, we do not assume symmetry). We denote
\[ \bm \pi = (\bm z, \bm \ell, \bm \sigma) \]
the collection of parameters.

Define the interaction strength (recall that $q$ is defined in \eqref{def:q} in link with $Q$)
\begin{equation} \label{edf:q*}
    q_* (\bm z) = \sum_{ j \ne k } q(|z_j-z_{k}|).
\end{equation}

We will always assume that the Lorentz boosts $\bs \ell$ are small and the relative distance of the center is large: that is, for some small $\delta_0 >0$,
\begin{equation}\label{lk:zk}
|\bm \ell| + q_*(\bm z) \le \delta_0.
\end{equation}

For $k=1, \dots, K$, define
\begin{equation}\label{def:gs}
Q_{k}(x)=\sigma_{k} Q(x-z_{k})\quad \text{and} \quad
\vec{Q}_{k}(x) =\begin{pmatrix} Q_{k}(x) \\ -({\ell}_{k} \cdot \nabla) Q_{k}(x) \end{pmatrix}.
\end{equation}
and similarly,
\begin{equation*}
Y_k(x)=\sigma_{k} Y(x-z_{k}),\quad
\vec Y_k^\pm(x)= \sigma_{k} \vec Y^\pm (x-z_{k}),\quad
\vec Z_k^\pm(x)= \sigma_{k} \vec Z^\pm (x-z_{k}).
\end{equation*}
Set
\begin{equation}\label{def:G}
S[\bm \pi] = \sum_{k=1}^K {Q}_k, \quad  \vec{S} [\bm \pi] = \sum_{k=1}^K \vec{Q}_k,
\quad 
G=f\left( \sum_{k=1}^K Q_{k}\right)- \sum_{k=1}^K f\left(Q_{k}\right).
\end{equation}
If the parameters $(\bm z, \bm \ell)$ are time dependent functions, we make the adapted modification to define $Q_k(t,x) = \sigma_k Q(x-z_k(t))$ etc.

The following lemma gathers some bounds on interactions of the solitons, mostly depending on $q_*(\bm z)$.

\begin{lemma}[Interactions]
Assume that \eqref{lk:zk} holds. Let $m >0$, and $1 \le k,k' \le K$ with $k \ne k'$.
\begin{enumerate}
\item\emph{Bounds.} For any $0<m'<m$,  then
\begin{gather}
\int |Q_{k} Q_{k'} |^m \lesssim e^{-m' |z_{k}-z_{k'}|},\\
\int \Big| F \left(\sum_{j=1}^K Q_j \right)-\sum_{j=1}^K F(Q_j)-\sum_{j \ne j'} f(Q_j)Q_{j'} \Big|
\lesssim q_*(\bm z)^{5/4} ,\label{tech3}
\end{gather}
		
\item \emph{Sharp bounds.}
\begin{align}
\int |Q_{k'}| |Q_k|^{1+m} &\lesssim q\left( |z_{k}-z_{k'} |\right)  \label{tech1}\\
|G| &\lesssim \sum_{j \ne j'} |Q_j|^{p-1}|Q_{j'}|,  \label{tech1.5} \\
\|G\|_{L^2} & \lesssim q_*(\bm z). \label{tech2}
\end{align}
  
\item\emph{Asymptotic.}
\begin{equation}\label{tech4}
\big|\langle f(Q_{k}),Q_{k'}\rangle-\sigma_{k}\sigma_{k'}c_{1} g_0	q(|z_{k}-z_{k'}|)\big|\lesssim |z_{k}-z_{k'}|^{-1}q(|z_{k}-z_{k'}|)
\end{equation}
where
\begin{equation}\label{def:kappa}
g_0:=\frac{1}{c_{1}}\int Q^{p}(x)e^{-x_1}\ud x>0\quad\mbox{and}\quad c_{1}:=\|\partial _{x_1}Q\|_{L^{2}}^{2}.
\end{equation}

\item\emph{Leading order interactions.}There exists a smooth function $g:[0,\infty)\to \m R$ such that, for any $0<\theta<\min(p-1,2)$ and $r>1$
\begin{equation}\label{on:gg}
|g(r)-g_0 q(r)|\lesssim r^{-1} q(r),\quad 
|g(r)-g(r')| \lesssim |r-r'|\max(q(r),q(r')),
\end{equation}
and
\begin{gather}\label{G:on}
\left|\langle G,\nabla Q_{k}\rangle - c_1\sum_{j \in \llbracket 1,  K \rrbracket \setminus \{  k \}}\sigma_{j}\sigma_{k} \frac{z_j-z_{k}}{|z_j-z_{k}|} g(|z_j-z_{k}|) \right| \lesssim q_*(\bm z)^\theta.
\end{gather}
\end{enumerate}
\end{lemma}

\begin{proof}
The proofs are a combination of that of \cite[Lemma 2.1]{CMYZ} and \cite[Lemma 3.2]{CMY}.

\medskip

\emph{Proof of (i)}
First we denote $z=z_k-z_{k'}$ 
By \eqref{Qdec},  $|Q(y)|+|\nabla Q(y)| \lesssim e^{-|y|}$, and thus
\[
\begin{aligned}
&\int\left|Q_k Q_{k'}\right|^m \mathrm{~d} y+\int\left(\left|\nabla Q_k\right|\left|\nabla Q_{k}\right|\right)^m \mathrm{~d} y \\& \lesssim \int e^{-m|y|} e^{-m^{\prime}|y+z|} \mathrm{d} y 
 \lesssim e^{-m^{\prime}|z|} \int e^{-\left(m-m^{\prime}\right)|y|} \mathrm{d} y \lesssim e^{-m^{\prime}|z|}.
\end{aligned}.
\]
\medskip

\emph{Proof of (ii)}
First we have the following point-wise estimate,
\begin{equation}\label{pp:Q}
Q(y-z) \lesssim(1+|y-z|)^{-\frac{d-1}{2}}  \mathrm{e}^{-|z|+|y|}.
\end{equation}
By using the \eqref{pp:Q} to deduce that
\begin{align*}
\MoveEqLeft \int |Q_{k'}| |Q_k|^{1+m} \ud x = \int Q(y-z) Q^{1+m}(y) \ud y\\
&\lesssim q(|z|)\int_{|y|<\frac 34 |z|} e^{|y|} Q^{1+m}(y) \ud y +e^{-|z|} \int_{|y|>\frac 34 |z|} e^{|y|} Q^{1+m}(y) \ud y\\
&\lesssim q(|z|).
	\end{align*}
	The second inequality also involves the basic fact that when $|y|<\frac 34 |z|$, $|y-z| \ge |z|-|y| \ge|z|-\frac{3}{4}|z|$ and the decay of $q$ (recall \eqref{Qdec}).

The fact \eqref{tech1.5} is derived by direct Taylor expansion, and together with \eqref{tech1}, it implies \eqref{tech2}.
\medskip

\emph{Proof of (iii)} We refer to \cite{CMYZ}[Lemma 2.1, (iii)].

\emph{Proof of (iv)} Using Taylor formula, there holds
\begin{align*}
\bigg|G-p|Q_{k}|^{p-1} \sum_{j \in \llbracket 1,  K \rrbracket \setminus \{  k \} }Q_{j} \bigg| & \lesssim |Q_{k}|^{p-2} \sum_{j \in \llbracket 1,  K \rrbracket \setminus \{  k \} }|Q_{j}|^{2} \\
& \qquad +
\sum_{j \in \llbracket 1,  K \rrbracket \setminus \{  k \} } \sum_{j' \in \llbracket 1,  K \rrbracket \setminus \{ j, k \} } |Q_{j'}|^{p-1}|Q_j| .
\end{align*}
Thus, using~\eqref{tech1}, for $i \in \llbracket 1, d \rrbracket$
\begin{align*}
\bigg|\langle G,\partial_{x_i} Q_{k}\rangle-\sum_{j \ne k}\sigma_{j}\sigma_{k} H_i(z_j-z_{k})\bigg|
&\lesssim \sum_{j'\neq j}\int |Q_{j}|^{p-1}|Q_{j'}|^{2} \ud x \lesssim q_*(\bm z)^\theta.
	\end{align*}
	where we denote 
	\begin{align}
        H_i (z) = \int  \partial_{x_i} (Q^p)(y)Q(y+z)dy,
        \end{align}
so that $\langle  \partial_{x_i} (Q_j^p),Q_{j'} \rangle = H_i(z_j - z_{j'})$.
As $Q$ has radial symmetry, there exists a function $g: [0,\infty) \rightarrow \mathbb{R}$ such that
\[ H(z)= c_1 \frac{z}{|z|} g(|z|), \]
and we refer the reader to \cite{CMYZ}[Lemma 2.1, (iv)] for details.
\end{proof}

\subsection{Decomposition close to the sum of symmetric solitary waves}

Let $G$ be a subgroup of $O_d(\m R)$ and $\tau: G \to \{ \pm 1\}$ be a homomorphism, and $\Omega$ be a finite set of cardinal $|\Omega| = K \ge 2$. We denote $\omega_1, \dots, \omega_K$ the elements of $\Omega$, and for $k = 1, \dots, K$, $\sigma_k = \sigma(\omega_k) \in \{ \pm 1 \}$ the signum. We assume that $G$ acts on $\Omega$, which we realize via the homomorphism $\chi : G \to \gh S_K$ defined by
\begin{equation} \label{def:chi}
\forall R \in G, \ \forall k = 1, \dots, K, \quad \omega_{\chi(R) k} = R \omega_k.
\end{equation}
Also, we assume that the sign assignment $\bs \sigma$ is $\tau$-equivariant (see \eqref{def:equiv_ts}), which writes
\begin{equation} \label{eq:tau_equiv}
\forall R \in G, \ \forall k = 1, \dots, K, \quad \sigma_{\chi(R) k} = \tau(R) \sigma_k.
\end{equation}
We recall that these assumptions are meant so that $S[\bs \omega, 0, \bs \pi] \in \q H_\tau$ has symmetry, see \eqref{eq:S_sym}.

We say that a family of $K$ points $\bm \xi = (\xi_k)_{1 \le k \le K} \in (\m R^d)^K$ is $G$-invariant if
\begin{equation}
    \label{def:G_chi_permu}
    \forall R \in G,\ \forall k=1, \dots, K, \quad R \xi_{k} = \xi_{\chi(R) k}.
\end{equation}

Denote 
\[ \q M = \{ \bm \xi = (\xi_k)_{1 \le k \le K} \in (\m R^d)^K : \bm \xi \text{ is } G\text{-invariant and for all } j \ne k, \ \xi_j \ne \xi_k \}. \]
Of course $\bs \omega \in \q M$. Observe that $\q M$ is locally a linear subspace: indeed, near a point $\bs \xi = (\xi_1, \dots, \xi_k) \in \q M$, a $K$-tuple in $\q M$ also has (pairwise) distinct elements, so that $\q M$ is defined locally by the linear constraints associated to \eqref{def:G_chi_permu}. In particular, $\q M$ is locally smooth (but the dimension can depend the connected component), and we denote as usual $T_{\bs \xi} M$ the tangent space at a point $\bs \xi \in \q M$ and $T\q M$ the tangent bundle.

\begin{lemma}[Modulation with symmetry]\label{lem:dec} 
There exists $\delta_0>0$ such that the following holds.

Let $\bs \xi \in \q M$ with $q_*(\bs \xi) \le \delta_0$. Denote $\bs \pi = (\bs \xi, 0, \bs \sigma)$ and
 \[ \q V(\bs \pi) = \left\{ \vec u = (u,v) \in \q H_\tau :   \| \vec u - \vec S[\bs \pi] \|_{H^1 \times L^2} \le \delta_0 \right\}. \]
Then there exist a neighbourhood $\q W(\bs \xi)$ of $(\bs \xi,0,0) \in \q M \times T_{\bs \xi} {\q M} \times \q H_\tau$ and a $\mc C^\infty$ map
\[ \Phi: \begin{array}[t]{rcl} \q V(\bm \pi) & \to & \q W(\bs \xi) \\
\vec u & \mapsto & (\bm z, \bm \ell, \vec \e), \quad \vec \e = (\e, \eta),
\end{array} 
\]
such that $(\bm z, \bm \ell, \vec \e)$ is the unique element of $\q W(\bs \xi)$ such that
\[ \vec u = \vec S[\bm z,\bm \ell, \bm \sigma] + \vec \e, \]
(in particular, $\bm z$ and $\bm \ell$ are $G$-invariant) and the following orthogonality conditions hold:
\[ \forall k=1, \dots, K,\  \forall i =1 , \dots, d,  \quad \langle \e, \partial_{x_i} Q_k \rangle = \langle \eta, \partial_{x_i} Q_k \rangle =0. \]
Furthermore,
\[ | \bm \xi - \bm z | +  |\bm \ell|  + \| \vec \e \|_{H^1 \times L^2} \lesssim   \| \vec u - \vec S[\bm \pi] \|_{H^1 \times L^2} + q_*(\bm \xi). \]
\end{lemma}

\begin{proof}
The principle of existence of $\Phi$ is classical, and a consequence of the implicit function theorem. First, let us give a description of the tangent space $T_{\bs \xi} {\q M}$ more amenable to computations. For this, decompose $\llbracket 1, K \rrbracket$ into orbits under $\chi$: $\llbracket 1, K \rrbracket = \bigsqcup_{i=1}^s \Omega_i$  ($s \le K$); up to reindexing, we can assume that $i \in \Omega_i$, and for each $i \in \llbracket 1, s\rrbracket$ and $j \in \Omega_i$, let $R_{ij} \in G$ be such that $\chi(R_{ij}) i = j$. Define
\[ E_i = \bigcap_{R \in \mathnormal{Stab} (i)} \ker (R - \Id) \quad \text{where} \quad \mathnormal{Stab} (i) = \{ R \in G : \chi(R) i = i \}. \]
Then $T_{\bs \xi} \q M \simeq \prod_{i=1}^s E_i$. Indeed  the constraints are equivalent to: for all $i \in \llbracket 1, s \rrbracket$, $\zeta_i \in E_i$ and $\zeta_j = R_{ij} \zeta_i$ for all $j \in \Omega_i$.

For each $i\in \llbracket 1, s \rrbracket$, let $(e_{ij})_j$ be a basis of $E_i$; consider the functional
\[ F: \begin{array}[t]{rcl}
\q H_\tau \times \q M \times T_{\bs \xi} \q M & \to & \prod_{i=1}^s \m R^{\dim E_i} \times \prod_{i=1}^s \m R^{\dim E_i} \\
(\vec u, \bm z, \bs \ell) & \mapsto & ((\langle \e, \partial_{e_{ij}} Q( \cdot - z_i) \rangle)_j, (\langle \eta, \partial_{e_{ij}} Q(\cdot - z_i) \rangle)_j)_{i \in \llbracket 1, s \rrbracket},
\end{array} \]
where $\vec \e = \vec u - \vec {S}[\bm z, \bm \ell, \bm \sigma]$. One computes that
$\partial_{(\bm z,\bm \ell)} F(\vec u,\bm z,\bm \ell)$ is invertible, with uniform bounds on the neighbourhood $\tilde{\q V} = \{ (\vec u, \bm z, \bm \ell) : \vec u \in \q V(\bm \pi), \ |\bm z - \bm \xi|, |\bm \ell| \le \delta_0 \}$. Indeed, the dimensions (equal to $2\dim T_{\bs \xi} \q M$) coincide, and on $\tilde {\q V}$ 
the (partial) Jacobian matrix 
\[ (\partial_{\bm z, e_{i,j}} F(\vec u,\bm z,\bm \ell), \partial_{\bs \ell,e_{i,j}} F(\vec u,\bm z,\bm \ell))_{i,j} \]
is a perturbation of uniform size $O(\| \vec \e \|_{H^1 \times L^2} + q_*(\bs \xi)+|\bm z - \bs \xi| + |\bs \ell|) = O(\delta_0)$, of 
\[ \begin{pmatrix} 
-\| \partial_{x_1} Q \|_{L^2}^2 I_{\dim T_{\bs \xi} \q M} & 0 \\
0 & \| \partial_{x_1} Q \|_{L^2}^2 I_{\dim T_{\bs \xi} \q M}
\end{pmatrix}. \]
Hence, the implicit function theorem applies on $\tilde {\q V}$, and there, $(\bm z, \bm \ell)$ satisfying $F(\vec u,\bm z,\bm \ell)=0$ can be expressed as a function of $\vec u$. From there, one gets $\vec \e \in \q H_\tau$.

We have obtained that
\begin{equation} \label{eq:Gsym_mod_Ei} 
\forall i=1,\dots, s,\  \forall v \in E_i, \quad \langle \e, v \cdot \nabla Q(\cdot - z_i) \rangle = 0,
\end{equation}
(and the same for $\eta$).
Let us now show that $\vec \e$ satisfies all the orthogonality conditions. Observe that for $R \in O_d(\m R)$ and $z,v \in \m R^d$, and as $Q$ has spherical symmetry,
\begin{multline} \label{eq:Gsym_mod}
\tau(R)  \langle \e, v \cdot \nabla Q(\cdot - z) \rangle = \langle \e(R^{-1} \cdot), v \cdot \nabla Q(\cdot - z) \rangle = \langle \e, v \cdot \nabla Q(R \cdot - z) \rangle \\
 = \langle \e, v \cdot R^{-1} \nabla [ Q(R (\cdot - R^{-1}z))] \rangle = \langle \e, R  v \cdot \nabla Q(\cdot - R^{-1}z) \rangle.
 \end{multline}
 
Let $i \in \llbracket 1, s \rrbracket$. Notice that if $R \in \mathnormal{Stab} (i)$, $\chi(R)i =i$ so that \eqref{eq:tau_equiv} implies $\tau(R) =1$. Then the equality \eqref{eq:Gsym_mod} gives
\begin{equation} \label{eq:Gsym_mod2}
\forall R \in \mathnormal{Stab} (i), \ \forall v \in \m R^d, \quad \langle \e, (R-\Id) v \cdot \nabla Q(\cdot - z_i) \rangle =0.
\end{equation}
Now, there exists a finite subset $(\tilde R_{ij})_j \subset \mathnormal{Stab} (i)$ such that 
\[ E_i = \bigcap_{j} \ker(\tilde R_{ij}^{-1} -  \Id). \]
Then
\[ E_i^\perp = \sum_j  (\ker(\tilde R_{ij}^{-1} -  \Id))^\perp = \sum_j  \Im(  \tilde R_{ij} - \Id). \]
On the one hand, consider $v \in E_i^\perp$. The above decomposition implies that there exists $v_1, \dots, v_j\in \m R^d$ such that $v = \sum_j (\tilde R_{ij} - \Id) v_j$.
But \eqref{eq:Gsym_mod2} implies that for any $j$, $ \langle \e,  (\tilde R_{ij} - \Id) v_j \cdot \nabla Q(\cdot - z_i)) \rangle =0$. Summing up in $j$, we infer that $ \langle \e, v \cdot \nabla Q(\cdot - z_i) \rangle =0$.

On the other hand, if $v \in E_i$, we recall \eqref{eq:Gsym_mod_Ei}. Decomposing over $\m R^d = E_i \oplus E_i^\perp$, we conclude that for any $v \in \m R^d$,
\[ \langle \e, v \cdot \nabla Q(\cdot - z_i) \rangle =0. \]
This completes the orthogonality conditions on $\e$. The same argument holds on $\eta$. Hence the orthogonality conditions are met. The estimates are a consequence of the implicit function theorem.
\end{proof}

\begin{lemma}[Modulation equations of a solution of \eqref{DLKG}] \label{le-dec}
Let $\delta_0 >0$ be as in Lemma \ref{lem:dec}, and let  $\vec u=(u,\partial_t u) \in \mc C([T_1,T_2], \q H_\tau)$ be a symmetric solution of~\eqref{DLKG} defined on an interval $[T_1 ,T_2 ]$
such that for all $t \in [T_1,T_2]$ there exists $\bm \xi(t) \in (\m R^d)^{K}$ with $q_*(\bm \xi(t)) < \delta_0$ and
\begin{equation}\label{for:dec}
\| \vec u(t) - \vec S[\bm \xi(t), 0, \bm \sigma] \|_{H^1 \times L^2} \le \delta_0.
\end{equation}
Then there exist unique $\mathscr C^1$ functions $(\bm z,\bm \ell): [T_1,T_2]  \to T\q M$, such that the solution $\vec{u}$ decomposes on $[T_1 ,T_2 ]$ as
\[ \vec u(t) = \vec S[\bm z(t),\bm \ell(t), \bm \sigma] + \vec \e(t), \]
and with the following properties on $[T_1 , T_2 ]$.
\begin{enumerate}
		\item \emph{Orthogonality and smallness.} For any $k =1, \dots, K $ and $i=1, \dots, d$,
		\begin{equation}\label{ortho}
		\langle \e, \partial_{x_i} Q_{k}\rangle=\langle \eta,\partial_{x_i} Q_{k}\rangle=0
		\end{equation}
		and
		\begin{equation}\label{eq:bound}
		\|\vve\|_\ENE + | \bm \ell| + q_*(\bm z) \lesssim\delta_0.
		\end{equation}
		
		\item \emph{Equation of $\vve$.} Let $\bm \pi = (\bm z, \bm \ell, \bm \sigma)$ then
		\begin{equation}\label{syst_e}\left\{\begin{aligned}
		\partial_t \e & = \eta + \md _{{\mathbf \e}}, \\
		\partial_t \eta &
		= \Delta \e-\e+f(S[\bm \pi]+\e)-f(S[\bm \pi])
		-2\alpha\eta + \md_{\eta}+G,
		\end{aligned}\right.\end{equation}
		where
		\begin{align*}
		\mathrm{Mod}_{\e} &=
		\sum_{k=1}^K \left(\dot z_k-
		{\ell}_{k} \right)\cdot\nabla Q_k,\\
		\mathrm{Mod}_{\eta} & = \sum_{k=1}^K \big( \dot{\ell}_{k} + 2\alpha {\ell}_{k}\big) \cdot \nabla Q_{k}-\sum_{k=1}^K  ({\ell}_{k}\cdot\nabla)(\dot{z}_{k}\cdot\nabla)Q_{k}.
		\end{align*}
		\item \emph{Equations of the geometric parameters.} For $k=1, \dots, K$,
		\begin{align}
		|\dot{z}_{k}-{\ell}_{k}|&\lesssim \|\vve \|^{2}_\ENE+ |\bm \ell|^2,
		\label{eq:z}\\
		|\dot{{\ell}}_{k}+2\alpha {\ell}_k|
		&\lesssim \|\vve\|_\ENE^{2} + |\bm \ell|^2 + q_*(\bm z).\label{eq:l}
		\end{align}
		\item{}\emph{Refined equation for $\ell_{k}$.}
		For any $1<\theta<\min(p-1,2)$, $k =1, \dots, K$
		\begin{equation}\label{eq:lbis}
		\Big| \dot{{\ell}}_{k}+2\alpha{\ell}_{k}+ \sum_{j \in \llbracket 1, K \rrbracket \setminus \{ k \}} \sigma_k \sigma_j  \frac{z_k - z_j}{|z_ k -z_j|}g(|z_k-z_j|)\Big|
		\lesssim \|\vve\|_\ENE^{2}+|\bm \ell|^{2}+  q_*(\bm z)^\theta.
		\end{equation}
		
		\item\emph{Equations of the exponential directions.} 
		For $k=1, \dots, K$, let
		\begin{equation}\label{def:a}
		a_k^{\pm} = \langle \vve,\vec Z_{k}^{\pm}\rangle.
		\end{equation}
		Then,
		\begin{equation}\label{eq:a}
		\left| \frac \ud{\ud t} a_k^{\pm}- \nu^\pm a_{k}^{\pm}\right|\lesssim \|\vve \|_\ENE^{2} + |\bm \ell|^{2} + q_*(\bm z).
		\end{equation}
	\end{enumerate}
\end{lemma}

\begin{proof}
(i) The decomposition for each $t \in [T_1,T_2]$ fixed follows immediately from the static decomposition Lemma \ref{lem:dec}; the regularity of $\bm z, \bm \ell$ is a consequence of smoothness the map $\Phi$ there, and of the regularity in $\vec u$.

\medskip

In the rest of this proof, we formally derive the equations of $\vve$ and the geometric parameters from the equation of $u$. The strategy is to write \eqref{ortho} as a non-autonomous ODE system and use Cauchy-Lipschitz theorem to conclude that the parameters are $\mc C^1$ functions of time in some interval.

\bigskip
	
(ii) First, by the definition of $\e$ and $\eta$,
		\begin{equation*}
		\partial_t\e =\partial_{t}u-\sum_{k=1}^K \partial_{t} Q_{k}
		=\eta+\sum_{k=1}^K \left(\dot{z}_{k}-\ell_{k} \right)\cdot\nabla Q_{k}.
		\end{equation*}
		Second,
		\begin{align*}
		\partial_{t}\eta&=\partial_{tt}u+\sum_{k=1}^K \partial_{t}\left(\ell_{k}\cdot\nabla Q_{k}\right)\\
		&=\Delta u-u+f(u)-2\alpha\partial_{t}u+\sum_{k=1}^K \dot{\ell}_{k} \cdot\nabla Q_{k}-\sum_{k=1}^K\left(\ell_{k} \cdot\nabla \right)\left(\dot{z}_{k}\cdot\nabla \right) Q_{k}.
		\end{align*}
		By~\eqref{eq:elliptic}, $\Delta Q_{k} -Q_{k}+f(Q_{k})=0$ and the definition of $G$, denoting $S = S[\bm \pi]$ for simplicity:
		\begin{align*}
		\Delta u-u+f(u)-2\alpha\partial_{t}u
		&=\Delta\e-\e+f (S+\e )-f (S ) -2\alpha\eta\\
		&\quad +2\alpha\sum_{k=1}^K \ell_{k}\cdot\nabla Q_{k}+G.
		\end{align*}
		Therefore,
		\begin{align*}
		\partial_t \eta & =
		\Delta \e-\e+ f(S+\e)-f(S)-2\alpha\eta \\
		&\quad +\sum_{k=1}^K\big(\dot{\ell}_{k}+2\alpha{\ell}_{k}\big)\cdot\nabla Q_{k}
		-\sum_{k=1}^K({\ell}_{k} \cdot\nabla)(\dot{z}_{k}\cdot\nabla)Q_{k} + G.
		\end{align*}

\bigskip

(iii)-(iv) We derive~\eqref{eq:z} from~\eqref{ortho}. For any $i \in \llbracket 1, d \rrbracket$ we have
		\begin{equation*}
		0=\frac \ud{\ud t}\langle\e,\partial_{x_i} Q_k \rangle=\langle\partial_t \e,\partial_{x_i} Q_k \rangle + \langle \e,\partial_{t} (\partial_{x_i} Q_k )\rangle.
		\end{equation*}
		Thus the scalar product of the $\e$ line of \eqref{syst_e} with $\partial_{x_i} Q_k$ gives
		\begin{align*}
		\langle\eta, \partial_{x_i} Q_k \rangle+\langle \md_{\e}, \partial_{x_i} Q_k \rangle
		-\langle\e, \dot{z}_k \cdot\nabla\partial_{x_i} Q_k \rangle=0.
		\end{align*}
		The first term is zero due to the orthogonality~\eqref{ortho}. Hence,
		\begin{align*}
		\MoveEqLeft (\dot{z}_{k,i}-\ell_{k,i})\|\partial_{x_i} Q\|_{L^2}^2 \\
  & = - \sum_{j \in \llbracket 1,K \rrbracket \setminus \{ k \}}  \int(\dot{z}_j - \ell_j)\cdot \nabla Q_j (x)\left( \partial_{x_i} Q_k(x)\right) \ud x
		+\left\langle\varepsilon, \dot{z}_k \cdot \nabla \partial_{x_i} Q_k \right\rangle.
		\end{align*}
		Thus, also using~\eqref{tech1} with $m=1$ and $m'=\frac 12$, we obtain
		\begin{align*}
		|\dot{z}_{k,i}-\ell_{k,i}|\lesssim \sum_{j \in \llbracket 1,K \rrbracket \setminus \{ k \}} |\dot{z}_j - \ell_j| e^{-\frac 12 |z|}+|\dot{z}_k - \ell_k| \|\e\|_{L^2}
		+|\ell_k| \|\e\|_{L^2}.
		\end{align*}
		Since $\|\vve \|_\ENE\lesssim \gamma$, this yields
		\begin{align*}
		|\dot{z}_{k,i}-\ell_{k,i}|\lesssim \sum_{j \in \llbracket 1,K \rrbracket \setminus \{ k \}} |\dot{z}_j - \ell_j| e^{-\frac 12 |z_{j} - z_{k})|}
		+ |\ell_k| \|\e\|_{L^2}.
		\end{align*}
		And thus, for small $\delta_0$, $ e^{-\frac 12 |z_{j} - z_{k})|} \lesssim q_*(\bm z) \lesssim \delta_0$, so that, summing up in $i \in \llbracket 1, d \rrbracket$ and $k \in \llbracket 1, K \rrbracket$, we get 
		\begin{equation*}
		\sum_{k=1}^K |\dot{z}_{k}-{\ell}_{k}| \lesssim |\bm \ell|  \|\vve \|_\ENE,
		\end{equation*}
		which implies~\eqref{eq:z}.

  \medskip
  
    Next, we derive~\eqref{eq:l}-\eqref{eq:lbis}. From~\eqref{ortho}, it holds
    \begin{equation*}
	0=\frac \ud{\ud t}\langle{\eta},\partial_{x_i} Q_k \rangle
	=\langle\partial_t \eta,\partial_{x_i} Q_k \rangle
		+\langle{\eta},\partial_{t} (\partial_{x_i} Q_k )\rangle.
    \end{equation*}
    Thus, by~\eqref{ortho} and~\eqref{syst_e}, we have
    \begin{align*}
	0 & =\langle\Delta\e-\e+f'(Q_k )\e,\partial_{x_i} Q_k \rangle +\langle f(S+\e)-f(S)-f'(S)\e,\partial_{x_i} Q_k \rangle\\
	&\quad +\langle (f'(S)-f'(Q_k))\e, \partial_{x_i} Q_k \rangle \\
	& \quad +\langle\md_{\eta}, \partial_{x_i} Q_k \rangle+\langle G,\partial_{x_i} Q_k \rangle -\langle\eta, (\dot{z}_k \cdot\nabla )\partial_{x_i} Q_k \rangle.
    \end{align*}
    Since $\partial_{x_i} Q_k$ satisfies \[ \Delta\partial_{x_i} Q_k-\partial_{x_i} Q_k + f'(Q_k )\partial_{x_i} Q_k=0, \]
    the first term is zero. Next, by Taylor expansion (as $f$ is $\mathcal C^2$), we have
		\[  f(S+\e)-f(S)-f'(S)\e = \e^2 \int_0^1 (1-\theta) f''(S+\theta \e) \ud \theta, \]
		and by the $H^1$ sub-criticality of the exponent $p>2$, we infer
		\begin{equation}\label{untard}
		\left|\langle f(S+\e)-f(S)-f'(S)\e,\partial_{x_i} Q_k \rangle \right|\lesssim \|\e\|_{H^1}^{2}.
		\end{equation}
		Then, again by Taylor expansion and $p> 2$,
		\begin{equation*}
		|f'(S)-f'(Q_k)| 
		\lesssim\left|Q_k\right|^{p-2} \sum_{j \in \llbracket 1, K \rrbracket \setminus \{ k \}} \left|Q_{j}\right| + \sum_{j \in \llbracket 1, K \rrbracket \setminus \{ k \}} \left| Q_{j} \right|^{p-1},
		\end{equation*}
  so that
  \[ |f'(S)-f'(Q_k)| |\partial_{x_i} Q_k| \lesssim q_*(\bm z) \sum_{j \in \llbracket 1, K \rrbracket \setminus \{ k \}} |Q_j|^{p-2}, \]
and so (as $p >2$), using Cauchy-Schwarz inequality,
	\begin{equation}\label{deuxtard}
		|\langle (f'(R)-f'(Q_k))\e, \partial_{x_i} Q_k \rangle|
		\lesssim q_*(\bm z) \|\e\|_{L^2} \lesssim \| \e \|_{L^2}^2 + q_*(\bm z)^2.
	\end{equation}
    Perform the scalar product of $\partial_{x_i} Q_k$ with the line for $\eta$ in \eqref{syst_e}: direct computations show that
    \begin{align*}
	\MoveEqLeft \left(\dot{\ell}_{k, i}+2 \alpha \ell_{k, i}\right) \left\|\partial_{x_i} Q\right\|_{L^2}^2 \\
        & = \sum_{j=1}^K \left\langle\left(\ell_j \cdot \nabla\right)\left(\dot{z}_j \cdot \nabla\right) Q_j, \partial_{x_i} Q_k\right\rangle -\sum_{j \in \llbracket 1, K \rrbracket \setminus \{ k \}} \int\left(\dot{\ell}_j+2 \alpha \ell_j\right) \cdot \nabla Q_j\left(\partial_{x_i} Q_j\right) \\
        & \quad -\left\langle f(S+\varepsilon)-f(S)-f^{\prime}(S) \varepsilon, \partial_{x_i} Q_k\right\rangle \\
		& \quad -\left\langle\left(f^{\prime}(S)-f^{\prime}\left(Q_k\right)\right) \varepsilon, \partial_{x_i} Q_k\right\rangle-\left\langle G, \partial_{x_i} Q_k\right\rangle+\left\langle\eta,\left(\dot{z}_k \cdot \nabla\right) \partial_{x_i} Q_k\right\rangle.
		\end{align*}
    By~\eqref{G:on}, we have, for any $1<\theta<\min(p-1,2)$,		
    \begin{equation*} 
        \left|\frac{\langle G,\partial_{x_i} Q_k \rangle}{\| \partial_{x_i} Q\|_{L^2}^{2}} -c_1\sum_{j \in \llbracket 1, K \rrbracket \setminus \{ k \}}   \sigma_k \sigma_{j} \frac{z_{k,i}-z_{j,i}}{|z_k-z_j|}g(|z_k-z _j|)
	\right| \lesssim q_*(\bm z)^\theta.
    \end{equation*} 
    Finally, by~\eqref{eq:z} again, we bound
		\begin{align*}
		|\langle\eta, (\dot{z}_k \cdot\nabla )\partial_{x_i} Q_k \rangle |
		\lesssim (|\dot{z}_k - \ell_k| + |\ell_k|) \|\vve\|_\ENE \lesssim \|\vve \|_\ENE^2+|\bm \ell|^2.
		\end{align*}
    Combining the above estimates, we have obtained
    \begin{multline*}
	\Big|\big(\dot\ell_k+2\alpha\ell_k\big)+\sum_{j \in \llbracket 1, K \rrbracket \setminus \{ k \}} \sigma_k \sigma_j \frac{z_k-z_j}{\left|z_k-z_j\right|} g\left(\left|z_k-z_j \right|\right)\Big| \\
        \lesssim \sum_{j \in \llbracket 1, K \rrbracket \setminus \{ k \}}  |\dot\ell_j+2\alpha \ell_j |e^{-m |z_k-z_j|} + |\bm \ell |^2+\|\vve\|_\ENE^2+  q_*(\bm z).
    \end{multline*}
		These estimates imply~\eqref{eq:lbis};~\eqref{eq:l} follows readily using~\eqref{on:gg}.
	
\bigskip
	
(v) By~\eqref{syst_e}, we have
		\begin{align*}
		\frac \ud{\ud t}a_k ^{\pm} & = \langle\partial_{t}\vve,\vec{Z}_k ^{\pm}\rangle
		+\langle\vve,\partial_{t}\vec{Z}_k^{\pm}\rangle\\
		&=(\zeta^\pm-2\alpha) \langle \eta , Y_k \rangle+\langle \Delta \e-\e+f'(Q_k) \e,Y_k\rangle
		\\
		&\quad + \langle f(S+\e)-f(S)-f'(S)\e, Y_k\rangle
		+ \langle (f'(S)-f'(Q_k))\e, Y_k \rangle\\
		&\quad +\langle G,Y_k \rangle+\zeta^\pm \langle \md_\e,Y_k \rangle
		+\langle \md_\eta,Y_k \rangle-\langle\vve, \dot{z}_k \cdot \nabla\vec{Z}_k ^{\pm}\rangle.
		\end{align*}
		Using ~\eqref{def:a}, and 
		\begin{align}
		\zeta^\pm-2\alpha&=\nu^\pm\\
		\LL Y&=-\nu_0^2 Y\\\nu_0^2 &= \nu^\pm \zeta^\pm,
		\end{align}
		 observe that
		\begin{equation*}
		(\zeta^\pm-2\alpha) \langle \eta , Y_k \rangle+\langle \Delta \e-\e+f'(Q_k) \e,Y_k \rangle
		= \nu^\pm a_k^\pm.
		\end{equation*}
		Using the decay properties of $Y$  and proceeding as before for~\eqref{untard} and~\eqref{deuxtard},
		\begin{multline*}
		|\langle f(S+\e)-f(S)-f'(S)\e, Y_k\rangle|+|\langle (f'(S)-f'(Q_k))\e, Y_k\rangle|\\
		\lesssim \| \e\|_{L^2}^2 +\sum_{j \in \llbracket 1, K \rrbracket \setminus \{ k \}} e^{-\sqrt{1+\nu_0^2}|z_k-z_j|} \lesssim \| \e\|_{L^2}^2 + q_*(\bm z).
		\end{multline*}
		Next, by \eqref{tech2}, $|\langle G,Y_k \rangle|\lesssim q_*(\bm z)$. 
    Last, by~\eqref{eq:z} and~\eqref{eq:l},
\begin{equation*}
    |\langle \md_\e,Y_k \rangle|+|\langle \md_\eta,Y_k \rangle|+|\langle\vve, \dot{z}_k \cdot \nabla\vec{Z}_k ^{\pm}\rangle| \lesssim \|\vve\|_\ENE^2 +  |\bm \ell|^2 +q_*(\bm z). \qedhere
\end{equation*}
\end{proof}

\subsection{Energy estimates}

For $\mu>0$ small to be chosen, we denote $\rho=2\alpha-\mu$. Consider the nonlinear energy functional
\begin{equation} \label{def:E}
\mathcal{E} =
\int \left( |\nabla\e|^2+ (1-\rho\mu )\e^{2}
+(\eta+\mu\e)^2- 2 [F (S +\e )-F (S )-f (S )\e] \right).
\end{equation}

\begin{lemma}\label{le:ener} In the context of Lemma \ref{le-dec}, there exists $\mu>0$ such that
	the following hold.
	\begin{enumerate}
		\item\emph{Coercivity and bound.}
		\begin{equation}\label{eq:coer}
		\mu \|\vve \|_\ENE^{2}-\frac{1}{2\mu} (|\bm a^+|^2 + |\bm a^-|^2)
		\leq \mathcal{E}\leq \frac 1\mu \|\vve \|_\ENE^{2}.
		\end{equation}
		\item\emph{Time variation.}
		\begin{equation}\label{eq:E}
		\frac \ud{\ud t}\mathcal{E}\le -2\mu\mathcal{E}
		+\frac{1}{\mu}\|\vve \|_\ENE \bigg[\|\vve \|_\ENE^{2}+ |\bm \ell|^{2}+ q_*(\bm z)\bigg].
		\end{equation}
	\end{enumerate} 
\end{lemma}

\begin{proof}
	See[\cite{CMYZ}, proof of Lemma 2.4]. 
\end{proof}

\begin{remark}
	The above lemma holds for any small enough $\mu >0$.  For further needs, we assume \[ \mu \le\min
	\left\lbrace 1,\alpha,|\nu|\right\rbrace 
	.\]
\end{remark}

\begin{lemma}[Expansion of the energy] \label{ene expand} 
In the context of Lemma \ref{le-dec}, we have 
\begin{multline}
\label{exp:energy}
 E(\vec{u})= K E(Q, 0) - c_1g_0 \sum_{j \ne k}  \sigma_j \sigma_k q(|z_j - z_k|) \\
 + O\left( \frac{q_*(\bm z)}{|\ln(q_*(\bm z))|} + \|\vec{\varepsilon}\|_{H^1 \times L^2}^2 + |\bm \ell|^2\right).
\end{multline} 
\end{lemma}

\begin{proof}
This is reminiscent of \cite[Lemma 3.12]{CMY}.
To start with, we expand $E\left(u, \partial_t u\right)$ by using the decomposition \eqref{for:dec} and integration by parts, the equation 
\[ -\Delta Q_k+Q_k-f\left(Q_k\right)=0, \] and the definition of $G$ and $S[\bm \pi]$ in \eqref{def:G}: this yields
\begin{multline*}
2E\left(u, \partial_t u\right) =2E(S[\bm \pi], 0) +\int\left|\partial_t u\right|^2-2 \int G \varepsilon \\
 +\int\left(|\nabla \varepsilon|^2+\varepsilon^2-2 F(S+\varepsilon)+2 F(S)+2 f(S) \varepsilon\right).
\end{multline*}

Then due to \eqref{tech2}, the sub-criticality of $p$ and Sobolev embedding, there holds
\[ 2 E\left(u, \partial_t u\right)=\int\left|\partial_t u\right|^2+2 E(S[\bm \pi], 0)+O\left(q_*(\bm z)\|\vec{\varepsilon}\|_{H^1 \times L^2}+\|\vec{\varepsilon}\|_{H^1 \times L^2}^2\right).
\]
Note that 
\[ \| \partial_t u \|_{L^2}= \| \eta-\sum_{k=1}^K\left(\ell_k \cdot \nabla\right) Q_k \|_{L^2} \lesssim \|\vec{\varepsilon}\|_{H^1 \times L^2} + |\bm \ell|. \]
Then using \eqref{tech3},\eqref{tech4} and  \eqref{on:gg}, we compute: 
\begin{align*}
E(S[\bm \pi], 0)= & KE(Q, 0) \\
& \quad +\sum_{j \ne k} \int\left[\nabla Q_j\cdot \nabla Q_{k}\right.  \left.+Q_j Q_{k}-f\left(Q_j \right) Q_{k}-f\left(Q_{k}\right) Q_j\right] \\ 
& \quad -\int\left(F(R)-\sum_{k=1}^K F\left(Q_k\right)-\sum_{j \ne k} f\left(Q_j\right) Q_{k}\right) \\ 
= & K E(Q, 0)-\sum_{j \ne k} \left\langle f \left(Q_j \right), Q_k \right\rangle+O(q_*(\bm z)^\frac{5}{4}) \\
=& K E(Q, 0)-c_1g_0\sum_{j \ne k} q(|z_j-z_k|) \\
& \quad +O\Big(\sum_{j \ne k}|z_j-z_{k}|^{-1}q(|z_j-z_{k}|\Big) 
\end{align*}
Notice that $r \gtrsim |\ln q(r)|$ for $r \ge 1$, so that the $O$ on the last line can be replaced by $q_*(\bm z)/|\ln(q_*(\bm z))|$  (which also bounds $q_*(\bm z)^2$), and we arrive at \eqref{exp:energy}.
\end{proof}

\subsection{Time evolution analysis}\label{S:2.4}

We introduce new parameters and functionals to analyse the time evolution of solutions in the framework of Lemma~\ref{le-dec}.
First, we set $\bm y = (y_k)_{1 \le k \le K} \in (\m R^d)^K$ with
\begin{equation*}
y_{k}=z_{k}+\frac{\ell_{k}}{2\alpha},
\end{equation*}
In view of the asymptotics \eqref{Qdec} and since
\[ |z_k-z_{j}|-|y_k-y_{j}|=O(|\bm \ell|)\ll 1, \] 
 we have,
\begin{align*}
\left| q_*(\bm z) - q_*(\bm y) \right| \lesssim |\bm \ell| q_*(\bm y). 
\end{align*}
In particular, $q_*(\bm z) \sim q_*(\bm y)$.
Then, we introduce notation for the damped components of the solution
\begin{equation*}
\mathcal F = \mathcal E + \BB,\quad \BB= |\bm \ell|^2 + \frac 1 {2 \mu}   |\bm a^-|^2 ,
\end{equation*}
and for all the components except distances
\begin{equation*}
\mathcal N = \bigg[ \|\vve\|_\ENE^2 + |\bm \ell|^2 \bigg]^{\frac 12}.
\end{equation*}
Finally, we define
\begin{equation*}
b = |\bm a^+|^2 = \sum_{k=1}^K (a^+_k)^2, \quad \mathcal M= \frac 1{\mu^2}\left(\mathcal F - \frac{b}{2\nu^+}\right).
\end{equation*}

We rewrite the estimates of Lemmas~\ref{le-dec} and~\ref{le:ener} with these notations, with errors in terms of $q_*(\bm y)$, $\q N$.

\begin{lemma}\label{le:new}
Let any $1<\theta<\min (p-1,\frac{5}{4})$.
In the context of Lemma~\ref{le-dec}, the following hold.
\begin{enumerate}
    \item \emph{Comparison with original variables.} For $j,k=1, \dots, K$,

    \begin{gather}
	\big| |y_k  - y_j| -|z_{k} - z_{j}| \big| \lesssim \mathcal N, \label{eq:new1} \\
	\mu \mathcal N^2 \le \mu \|\vve\|_\ENE^2 + |\bm \ell|^2 \le \mathcal F + \frac {b}{2\mu}
	\lesssim\mathcal N^2.\label{eq:new2}
    \end{gather}
		
    \item\emph{ODE for the distances between solitary waves.} 
    For any $k = 1, \dots, K$, the equation for the evolution of $y_k$ is
    \begin{equation} \label{eq:dist}
	\dot{y}_{k} = - \frac{1}{2\alpha} \sum_{j \in \llbracket 1, K \rrbracket \setminus \{ k \}} \sigma_j \sigma_k \frac{y_j -y_k}{|y_j-y_k|}g(|y_j -y_k|)+O\big(\mathcal{N}^{2}+q_{*}(\bm y)^{\theta}\big),
    \end{equation}
    where $\theta$ is as in \eqref{eq:lbis}.
		
    \item\emph{Exponential instability.}
    \begin{equation}\label{eq:b}
        | \dot b - 2 \nu^+ b |\lesssim \mathcal N^3+ \mathcal N q _*(\bm y).
    \end{equation}

    \item\emph{Damped components.}
    \begin{align}\label{eq:damped1}
        \frac \ud{\ud t}\mathcal F + 2 \mu \mathcal F & \lesssim \mathcal N^3+ \mathcal N  q _*(\bm y),	\\
        \label{eq:damped2} 
        \frac \ud{\ud t}\BB + 2 \mu \BB & \lesssim \mathcal N^3+ \mathcal N q _*(\bm y).
   \end{align}
    
    \item\emph{Liapunov type functional.} 
    \begin{equation} \label{eq:Mbis}
        \frac \ud{\ud t} \mathcal M +\mathcal N^2 \lesssim q_*(\bm y)^2.
    \end{equation}
\end{enumerate}
\end{lemma}
\begin{proof}
	\emph{Proof of~\eqref{eq:new1}.} It follows
	from the triangle inequality that
	\begin{equation*}
	\big||y_{k} - y_{j}| - |z_{k} - z_{j}| \big| \le |y_k - z_{k} - (y_j-z_j)|\le \frac{|\ell_k- \ell_j|}{2\alpha }\lesssim \mathcal N.
	\end{equation*}
	The second part of~\eqref{eq:new1} then follows from~\eqref{on:gg}.
	
		\emph{Proof of~\eqref{eq:new2}. } It follows readily from~\eqref{eq:coer}.
	
	\emph{Proof of~\eqref{eq:dist}.} First, from~\eqref{eq:z},~\eqref{eq:lbis} and~\eqref{eq:new1}, we note ($\theta$ is as in \eqref{eq:lbis})
	\begin{equation}\label{for:y}\begin{aligned}
	\dot y_k &= \dot z_k + \frac{\dot \ell_k}{2\alpha}={\ell_k}+\frac{\dot{\ell_k}}{2\alpha}+O(\mathcal{N}^2)
	\\&=-\sum_{j \in \llbracket 1, K \rrbracket \setminus \{ k \}} \frac{\sigma_k \sigma_j}{2\alpha} \frac{z_k-z_j}{|z_k-z_j|} g(|z_k-z_j|)+O(\mathcal N^2+ q_*(\bm y)^\theta)
	\end{aligned}
	\end{equation}
 
    \emph{Proof of~\eqref{eq:b}. } 
    There hold
	\begin{align*}
		 \dot b - 2 \nu^+ b =2 \bm a^+ \cdot \left( \frac \ud{\ud t} \bm a^{+}- \nu^+ \bm a^{+}\right) 
	\end{align*}
    Then \eqref{eq:b} follows from~\eqref{eq:a} and 
     \[ |\bm a^+|\lesssim \|\vve\|_\ENE \le \mathcal N. \]

    \emph{Proof of \eqref{eq:damped1}-\eqref{eq:damped2}. }
    From the expression of $\mathcal F$ and then~\eqref{eq:E},~\eqref{eq:l} and~\eqref{eq:a}
    \begin{align*}
	\frac \ud{\ud t} \mathcal F
	& = \frac \ud {\ud t} \mathcal E + 2 \bm \ell \cdot \bm \ell + \frac 1 \mu \bm a^- \cdot \bm a^- \\
	& \le -2\mu \mathcal E -4\alpha |\bm \ell|^2 + \frac{\nu^-} \mu |\bm a^-|^2
	+ O(\q N^3 +  q_*(\bm y) \q N).
    \end{align*}
    Since $0<\mu<\alpha$ and $0<\mu<|\nu^-|$ we obtain~\eqref{eq:damped1} for $\mathcal F$.
    
    The  proof of $\BB$ follows the same procedure.

    \emph{ Proof of \eqref{eq:Mbis}.}
    First, it follows from combining~\eqref{eq:b} and~\eqref{eq:damped1} that
	\begin{equation*}
	\mu^2 \frac \ud{\ud t} \mathcal M \leq -2\mu \left(\mathcal F + \frac b{2\mu}\right) + O(\mathcal N^3+ q_*(\bm y) \mathcal N),
	\end{equation*}
	Now, from~\eqref{eq:new2},
	we observe that
	\begin{equation*}
	\mathcal N^2
	\le \frac 1{\mu} \left(\mathcal F+\frac{b}{2\mu}\right)
	\le -\frac 1{2} \frac \ud{\ud t} \mathcal M + O(\mathcal N^{3}+ q_*(\bm y) \mathcal N),
	\end{equation*}
	and from there, for some $C>0$,
	\[ 
	\frac \ud{\ud t} \mathcal{M} + 2 \q N^2 \le C \q N^3+ C q_*(\bm y) \q N.
	\]
	As $q_*(\bm y) \q N  \le 1/(2C) \q N^2 + 2C q_*(\bm y)^2$, taking $\q N$ small enough gives \eqref{eq:Mbis}.
\end{proof}

\section{Construction of multi-solitons} \label{sec:const}

Denote $\ds \beta =\langle Y, Z \rangle^{-1} > 0$.

\begin{lemma}[Modulation at initial time]\label{le-W} 
In the context of Lemma \ref{lem:dec}, for $\bm \pi = (\bm z, \bm \ell, \bm \sigma)$ such that
\[ q_*(\bm z) \le \delta_0 \]
is small enough, there exist linear maps
	\begin{equation} \label{est:BV}
	 B= B_{\bm \pi}: \m R^{K} \rightarrow \m R^K \quad \text{and} \quad V = V_{\bm \pi}: \m R^{K} \to (\m R^d)^K,
	\end{equation}
	depending smoothly on $\bm \pi$ and satisfying
\[
	\| B-\beta \operatorname{Id}\| + \| V \| \lesssim q_*(\bm z)^{1/2}, \]
	and such that the function $W \left(\bm a \right): \mathbb{R}^d \rightarrow \mathbb{R}$ defined for $\bm a \in \m R^K$ by
\[
	W \left(\textbf{a}\right) :=  \sum_{k=1}^K  \left[ B (\bm a )_k Y_k + V\left(\bm a\right)_k \cdot \nabla Q_k \right]
	\]
	satisfies the orthogonality conditions
	\[ \forall k=1, \dots, K, \forall i = 1, \dots, d,  \quad 
	\left\langle W \left(\textbf{a}\right), \partial_{x_i} Q_k \right\rangle=0, \quad\left\langle W \left( \bm a \right), Y_k \right\rangle=\beta a_k ,
	\]
	In particular, setting
\[ \vec{W} (\textbf{a}) = \begin{pmatrix}
	W (\textbf{a}) \\
	\nu^{+} W (\textbf{a})
	\end{pmatrix}, \quad \text{it holds for all } k=1, \dots, K,   \quad \big\langle \vec{W}_{\bm \pi}\left(\textbf{a}\right), \vec{Z}_k^{+} \big\rangle = a_k. \]
 Furthermore, if $\bm a$ has $\tau$-symmetry in the sense that $a_{\chi(R) k} =  \tau(R) a_{k}$ for all $k=1, \dots, K$ and $R \in G$ (cf. \eqref{eq:tau_equiv}), then $\vec W(\bm a) \in \q H_\tau$.
\end{lemma}

\begin{proof} 
    Let $\bm a \in \m R^K$ and define for $\bm b \in \m R^K$ and $\bm v \in (\m R^d)^K$
    \begin{equation}\label{def-W}
        \tilde W\left(\bm b, \bm v \right)(x):=\sum_{k=1}^K \sigma_k \bigg[ b_k Y_k(x)+ \sum_{i=1}^d v_{i, k} \partial_{x_i} Q_k(x)\bigg].
    \end{equation}
    Consider the linear system of $K \cdot (d+1)$ equations
    \begin{equation} \label{eq:lin_sys_init}
         \langle \tilde W(\bm b, \bm v), Y_k \rangle=\beta a_k, \quad
         \langle \tilde W(\bm b, \bm v), \partial_{x_i} Q_k \rangle =0,
    \end{equation} 
    for $i = 1, \dots, d$ and $k=1, \dots,K$. Our goal is to solve for $\bm b, \bm v$ in function of $\bm{a}$.
		
    Observe that there hold, for any $i \ne i'$, the relations 
    \[ \langle Y, Z \rangle =\beta, \quad \langle Y,\partial_{x_i} Q\rangle =0, \|\partial_{x_i} Q \|_{L^2} = \frac{1}{d} \| \nabla Q \|_{L^2}^2, \quad \text{and} \quad \langle \partial_{x_{i'}} Q,\partial_{x_i} Q \rangle =0, \]
    and for any $i,i$ and $j \ne k$, the estimates 
    \[ \langle \partial_{x_{i'}} Q_j,\partial_{x_i} Q_{k}\rangle, \ \langle Y_k,\partial_{x_i} Q_{j}\rangle,\ \langle Y_k, Y_j \rangle  = O(q(|z_k-z_{j}|)^{1/2}),  \]
    (see~\eqref{tech1}). Therefore, the system \eqref{def-W}-\eqref{eq:lin_sys_init} writes
    \[
        M \begin{pmatrix} 
        \bm b \\
        \bm v
        \end{pmatrix} = \begin{pmatrix} 
        \bm a \\
        0
        \end{pmatrix}, \quad \text{where} \quad M = \begin{pmatrix} \beta^{-1} I_{K} & 0 \\
        0 & \frac{1}{d} \| \nabla Q \|_{L^2}^2 I_{d K}
    \end{pmatrix} + O(q_*(\bm z)^{1/2}). \]
    $M$ is invertible as soon as $q_*(\bm z)$ is small enough, with
    \[ M^{-1} = \begin{pmatrix} \beta I_{K} & 0 \\
        0 & d \| \nabla Q \|_{L^2}^{-2} I_{d K} 
    \end{pmatrix} + O(q_*(\bm z)^{1/2}). \]
    Then $W(\bm a) = \tilde W(M^{-1}(\bm a,0))$ is the unique solution, and the matrices $B$ and $V$ are defined as the first $K$ columns of $M^{-1}$:
    \[ M^{-1} = \begin{pmatrix} B & * \\
        V & *
    \end{pmatrix}. \]
    This also proves the estimates \eqref{est:BV}.

    Finally, for the case with symmetry, observe that if $\bm a$ with $\tau$-symmetry is completely defined by the value of its entries $a_k$ for indices $k$ running on a system of representants of the orbits of $\llbracket 1, K \rrbracket$ under $G$ (via $\chi$). Therefore, one can repeat the proof, restricting to these indices (for $\bm a$, $\bm b$ and $\bm v$): this gives existence of a $\tau$-symmetric $\bm b$, $\bm v$, and one clearly has $\vec{\tilde W}(\bm b, \bm v) \in \q H_\tau$.
\end{proof}

\begin{proposition}\label{pr:exist}
    Let $0 \le \delta < \delta_0$ be small enough. In the context of Theorem \ref{thm:abs}, for $\Omega$, $G$, $\tau$ and $\bm \sigma$, assume that:
    \begin{equation}\label{eq:pr:10}
	\begin{cases}
	q_*(\Omega) \le \delta,\\
	\bm \ell(0) = (\ell_k(0))_{1 \le k \le K} \in \mathcal B_{(\RR^{d})^{K}}(\delta) \ G \text{-invariant},\\
	\vvep(0)\in \mathcal B_{\q H_\tau}(\delta) \text{ satisfying } \langle \vvep(0),\vec Z_k^+(0)\rangle=0 \ \text{ for all } k =1, \dots, K. 
	\end{cases}
    \end{equation}
    Then there exists $\bm {\gh a}^+ = (\gh a^+_{k})_{1\le k \le K} \in  \bar{\mathcal B}_{\RR^{K}}(\delta^\frac 54)$ such that the solution $\vec u$  of~\eqref{DLKG} with initial data 
	\begin{equation*}
	\vec u(0)=  \vec S[ \bs \omega, \bm \ell(0), \bm \sigma] + \vec W(\bm{\gh a}^+) + \vvep(0)
	\end{equation*}
	is defined globally for positive times and satisfies
\begin{equation} \label{conv:prop} 
\forall t \ge 0, \quad \vec u(t) = \vec S[\bm z(t),0,\bm \sigma] + \vec \e(t), \quad \text{with} \quad \| \vec \e(t) \|_{\q H_\tau} \lesssim t^{-1},
\end{equation}
and for some point $\tau_\infty \in \m R^d$ and some rotation $R_\infty \in O_d(\m R)$, we have for all $k =1, \dots, K$:
\begin{align}  \label{conv:prop2}
 z_k(t) = \lambda_{\Omega} r(t) R_\infty \omega_k  + \tau_\infty \quad \text{where }   \quad r(t) = \ln t - \frac{d-1}{2} \ln \ln t + c_0 + O \left( \frac{\ln \ln t}{\ln t} \right).
 \end{align}
\end{proposition}

\begin{remark}
$c_0 =\ln \left( \frac{\kappa \gamma g_0}{2\alpha} \right)$ depends on the parameters of the PDE ($d,p,\alpha$), and on $\Omega$ but solely through the constant $\gamma$ appearing in \eqref{symmetry}.

$\lambda_\Omega$ is the inverse of the minimal distance between two vortices (see \eqref{def:lambda_O}): it depends only on the configuration $\Omega$. At least up to the first 3 main orders, the dependence to $\Omega$ is only through $\lambda_\Omega$ and $\gamma$.
\end{remark}

\begin{proof}

\emph{Step 1. Bootstrap for a given initial data}

\bigskip

Given $\bm {\gh a}^+ = (\gh a_k^+)_{1 \le k \le K} \in \q B_{\m R^{K}}(\delta^{5/4})$, we consider the solution $\vec u(t)$ of \eqref{DLKG} with initial data
\begin{align} \label{def:ini}
\vec u(0) = \vec S[\bm z(0), \bm \ell(0), \bm \sigma] + \vec W(\bm {\gh a}^+) + \vvep(0).
\end{align}

\emph{Decomposition.}
For any $T > 0$ such that the solution $\vec u$ satisfies~\eqref{for:dec} on $[0,T]$, we decompose it according to Lemma~\ref{le-dec}: for $t \in [0,T]$ we therefore dispose of $\bm z(t), \bm \ell(t), \vec \e(t)$, and from there, $\bm y(t)$ and $\bm a^\pm(t)$.

As $\vec u$ is $G$ invariant, $\bm z$ and $\bm \ell$ are $G$ invariant too, and $\bm y$ 
as well (for each time $t \in [0,T]$).

Also, due to the properties of the function $W$ in Lemma~\ref{le-W} and the orthogonality properties \eqref{ortho} of  $\vvep(0)$ assumed in~\eqref{eq:pr:10}, the initial data $\vec u(0)$ is modulated, in the sense that $ (\bm{z}(0),\bm{\ell}(0))$ and
\begin{equation*}
    \vve(0)=\vec W(\bm {\gh a}^+) + \vvep(0),
\end{equation*}
are the parameters of the decomposition of $\vec u(0)$. Due to Lemma~\ref{le-W}, for $k =1, \dots, K$,
	\begin{equation*}
a_k^+(0) =  \langle \vve(0),\vec Z_{k}^+(0)\rangle =  \langle \vec W(\bm {\gh a}^+), Z_k^+ \rangle =  \gh a_k^+,
\end{equation*}
Therefore, it holds from~\eqref{eq:pr:10}
\begin{equation}\label{at:0}
    \mathcal N(0)\lesssim \delta,\quad q_*(\bm z(0))\lesssim \delta^2 \quad \text{and} \quad b(0) \le \delta^{5/2}.
\end{equation}

\bigskip
 
Due to the rigidity of $G$, there exist $R_\sharp \in \mc C^1([0,T], O_d(\m R))$, $\lambda_\sharp \in \mc C^1([0,T], \left]0,+\infty \right[)$ and $\tau_\sharp \in \mc C^1([0,T],\m R^d)$ such that
\begin{equation} 
\forall j=1, \dots, K, \quad y_j(t) = \lambda_\sharp (t) R_\sharp(t) \omega_j + \tau_\sharp(t).
\end{equation}
As a consequence, for all $j=1, \dots, K$ and $t \in [0,T]$
\begin{equation} \label{eq:rigid}
|y_{j}(t) - y_k(t)| = \lambda_\sharp(t) |R_\sharp(t) (\omega_j - \omega_k)| = \lambda_\sharp(t) |\omega_j - \omega_k|.
\end{equation}

Consider the minimal distance between two (distinct) points of $\Omega$: it is realized by $\omega_{k_*}$ and $\omega_{k_\dagger}$ (with $k_*, k_\dagger \in \llbracket 1, K \rrbracket$ and $k_* \ne k_\dagger$):
\[ \forall j,k \in  \llbracket 1, K \rrbracket, j \ne k , \quad |\omega_j - \omega_{k_*}| \ge |\omega_{k_*} - \omega_{k_\dagger}|. \]
Now let $\Omega_* = \{ k \in  \llbracket 1, K \rrbracket : |\omega_k - \omega_{k_*}| = |\omega_{k_\dagger} - \omega_{k_*}| \}$ be the set of indices of nearest neighbours to $\omega_{k_*}$ and let
\begin{equation} \label{def:lambda_O} 
 \lambda_\Omega := \frac{1}{|\omega_{k_*} - \omega_{k_\dagger}|} \quad \text{and} \quad \vartheta := \min \left\{ \frac{|\omega_j  - \omega_{k_*}|}{|\omega_{k_*} - \omega_{k_\dagger}|} : j \in \llbracket 1,K \rrbracket \setminus \Omega_* \right\} > 1 \end{equation}
be respectively the inverse of the minimal distance and the minimal ratio of the distance for a non-neighbouring point to $\omega_{k_*}$. For later purposes we assume 
\[ 1 < \theta \le \vartheta. \]

Then for any $k \in \Omega_*$
\[ \forall t \in [0,T], \quad  | y_k(t) - y_{k_*}(t)| = r(t) \quad \text{where} \quad r(t) = \frac{\lambda_\sharp(t)}{\lambda_\Omega}. \]
And for any $j \notin \Omega_*$, 
\[ |y_{j}(t) - y_{k_*}(t)| = \lambda_\sharp(t) |\omega_j - \omega_{k_*}|  \ge \vartheta r(t). \]
In particular, we infer that $q_*(\bm y(t)) = q(r(t))$ and
\begin{align} \label{est:non_nearest_neighbour} 
\forall j \notin \Omega_*, \quad q(|y_{j}(t) - y_{k_*}(t)|) \lesssim q_*(\bm y(t))^\vartheta \lesssim  q_*(\bm y(t))^\theta.
\end{align}

As $\Omega$ has center of mass $0$, notice that the center of mass of $\bm y(t)$ is
\[ \frac{1}{K} \sum_{k=1}^K y_k(t) = \lambda_\sharp(t) R_\sharp(t) \sum_{k=1}^K \omega_k + \tau_\sharp(t) = \tau_\sharp(t). \]
It follows that for all $k=1, \dots, K$,
\[ y_k(t) -\tau_\sharp(t) = \lambda_\sharp(t) R_\sharp(t)\omega_k  = \lambda_\Omega r(t) R_\sharp(t) \omega_k. \]

\bigskip

\emph{Bootstrap estimates.}
Let $C_* >0$ be a large constant to be fixed later. We introduce the following bootstrap estimates
\begin{equation}\label{BS-2}
    \mathcal N(t)  \le \frac{C_*}{t+\delta^{-1}},  \qquad q(r (t))  \le \frac{C_*}{t+\delta^{-2}} \qquad \text{and} \quad
    b(t)  \le  \frac{1}{(t+\delta^{-1})^2}. 
	\end{equation}
	We set
	\begin{equation*}
	\TS=\sup\left\{ t\in [0,\infty)\hbox{ such that \eqref{BS-2} holds on $[0,t]$}\right\}\geq 0.
	\end{equation*}
    In the bootstrap argument below, $C_* \ge 1$ will be chosen large relative to all the implicit constants appearing in the inequalities. Then $\delta>0$ will be any sufficiently small number, chosen so small that in particular $C_* \delta \le 1$. We therefore emphasize that in the computations below, the implicit constant does not depend on $C_*$ nor $\delta$.
	
	\bigskip
	
	\emph{Estimates on the damped components.}
	The estimate on $\mathcal N$ is strictly improved on $[0,\TS]$  in the following way.\\
	From~\eqref{eq:damped1} and, it holds on $[\TD,\TS)$,
	\begin{align*}
	\frac{d}{dt}\left[ e^{2\mu t}\mathcal F\right] &=e^{2\mu t}\left(\frac{d}{dt}\mathcal F+2\mu \mathcal F \right)\\
	&\lesssim e^{2 \mu t} \left( \mathcal N^3+\mathcal N q_*\right) \\&\lesssim e^{2\mu t}\left(C_*^3  \delta^3 (1+\delta t)^{-3}+C_*C_q\delta^{\frac 52} (1+\delta t)^{-2}\right) \\&\lesssim e^{2\mu t}\delta^{2} (1+\delta t)^{-2},
	\end{align*}
	for $\delta$ small enough. Thus, integrating the above estimate on $[0,t]$,
	for any $t\in [0,\TS)$, 
	\begin{equation*}
	\mathcal F(t)-e^{-2\mu t}\mathcal F(0) \lesssim \delta^{2}\int_{0}^t e^{-2\mu(t-s)} (1+\delta s)^{-2} \ud s \lesssim  \delta^{2}  (1+\delta t)^{-2}
	\end{equation*}
	From~\eqref{eq:new2} and \eqref{BS-2}, 
	\[ \mathcal F(t)\lesssim  \delta^{2}  (1+\delta t)^{-2}+e^{-2\mu t}\mathcal F(0). \]
	Now by \eqref{BS-2} with the fact that $e^{-2\mu t} \le (1+\mu t)^{-2}$, then we have
	\[e^{-2\mu t}\mathcal F(0) \lesssim e^{-2\mu t}\delta^2\lesssim  \delta^2(1+\delta t)^{-2}. \]
	In particular, by \eqref{eq:new2} and still denote the universal constant as $C_{*}$ we obtain
	\begin{equation*}\label{bnf2}
	\|\vve\|_\ENE^2 \lesssim \mathcal F+b
	\lesssim C_* \delta^{2} (1+\delta t)^{-2}
	\end{equation*}
	Arguing similarly for the quantity $\BB$, we have
	\begin{equation*}
	\sum_{k=0}^{K} |\ell_\omega|^2 + \sum_{k=0}^{K} (a_k^-)^2\lesssim \BB \lesssim C_* \delta^{2} (1+\delta t)^{-2}.
	\end{equation*}
	Hence we obtain, for all $t \in [\TD, \TS)$,
	\begin{equation}\label{est:N}
	\mathcal N(t) \lesssim \sqrt{C_*}\delta  (1+\delta t)^{-1} \lesssim \frac{\sqrt{C_*}}{t+\delta^{-1}}.
	\end{equation}
	For $C_*$ large enough, this strictly improves the bootstrap estimate \eqref{BS-2} on~$\mathcal N$ on the interval $[0,\TS)$.
	
        \bigskip
	
	\emph{Estimate on the distance.}
	


Let us first derive the motion for the center of mass, which we recall is $\frac{1}{K} \sum_{k=1}^K y_k = \tau_\sharp$. For this recall \eqref{eq:dist}: for any $k=1, \dots, K$,
\[ \dot{y}_{k}(t)  = \sigma_{k} \sum_{j=1}^K \sigma_j \frac{y_j - y_{k}}{|y_j - y_{k}|} g(|y_j - y_{k}|) +   O(\mathcal{N}^2+ q(r(t))^\theta)
\]
Averaging in $k$, we see that the sum term cancels out, and so
\begin{equation} \label{est:translation}
\dot{\tau_\sharp}(t) = O(\mathcal{N}^2+ q(r(t))^\theta).
\end{equation}
Using again \eqref{eq:dist} on $k =1, \dots, K$, together with \eqref{est:non_nearest_neighbour} and  \eqref{symmetry} and denoting $\Omega_k$ the set of indices of nearest neighbours of $\omega_k$, we have
\begin{align*}\label{equ:dyk}
	\dot{y}_k(t) - \dot \tau_\sharp(t)& = \frac{\sigma_{k}}{2\alpha}  \sum_{j \in \llbracket 1, K \rrbracket \setminus \{ k \}} \sigma_j \frac{y_j - y_k}{|y_j - y_k|} g(|y_j - y_k|) +   O(\mathcal{N}^2+ q_*(r(t))^\theta) \\
	& =  \frac{\sigma_{k}}{2\alpha} \sum_{j \in \Omega_k} \sigma_j R_\sharp(t) \frac{\omega_j - \omega_k}{|\omega_j - \omega_k|} g(r(t)) +   O(\mathcal{N}^2+ q_*(r(t))^\theta) \\
 	& =  \frac{\gamma \lambda_\Omega}{2\alpha} g(r(t)) R_\sharp(t) \omega_k  +  O(\mathcal{N}^2+ q_*(r(t))^\theta)
\end{align*}
Now, as $y_k(t) - \tau_\sharp(t) = \lambda_\Omega r(t) R_\sharp(t) \omega_k$, we arrive at
\begin{equation} \label{est:rot}
\dot r R_\sharp \omega_k + r \dot R_\sharp \omega_k = \frac{\gamma}{2 \alpha} g(r(t)) R_\sharp \omega_k + O(\mathcal{N}^2+ q_*(r(t))^\theta).
\end{equation}
Taking the scalar product with $R_\sharp (t) \omega_k$, we infer the dynamical equation of $r$:
 \begin{align}
	\dot r(t) & =  \frac{\gamma}{2 \alpha} g(r(t)) + O(\mathcal{N}^2+ q_*(r(t))^\theta) \label{eq:drk2} \\
        & =\frac{\gamma g_0}{2\alpha} q(r(t)) + O(\mathcal{N}^2+ r(t)^{-1} q(r(t))).
\end{align}
(the last equality comes from \eqref{on:gg}). 

Then using the fact that
	$
	|q(r)+q'(r)|\lesssim r^{-1}q(r),
	$
	we can integrate this differential inequality: indeed there hold
	\begin{equation}
	\frac{d}{dt}\left(\frac{1}{q(r(t))} \right) =\frac{\gamma  g_0}{2\alpha} +O\left( \frac{\mathcal{N}^2}{q(r(t))}+r(t)^{-1}\right).
	\end{equation}
	which means that for some $C >0$
 \begin{equation}\label{equ:dis-sym}
	\left| \frac{d}{dt}\left(\frac{1}{q(r(t))} \right) -\frac{\gamma g_0}{2\alpha}\right|  \lesssim \left( \frac{\mathcal{N}^2}{q(r(t))}+r(t)^{-1} \right) .
	\end{equation}
	We consider the refined distance 
	\[
	D(t)=\frac{1}{q(r(t))}\cdot e^{-C\mathcal{M}(t)},
	\]
	then we have,
	\begin{align*}
	e^{C\mathcal{M}(t)} \dot D(t) &=\frac{d}{dt} \left[ \frac{1}{q(r(t))}\right] -\frac{{C}}{q({r}(t))} \dot{\mathcal M} \\
         & \ge \frac{\gamma g_0}{2\alpha}-C\frac{\mathcal{N}^2}{q(r(t))}+C\frac{\mathcal{N}^2}{q(r(t))}-C {r(t)}^{-1}-{C}q^2({r(t)}).
	\end{align*}
	Note that $r(t) \gtrsim |\ln \delta|$ and thus,
	\[
	\dot D(t) \geqslant \frac{\gamma g_0}{2\alpha} -\frac{{C}}{|\ln \delta|} -{C} \delta^2 \geqslant \frac{\gamma g_0}{4 \alpha},
	\]
    if $\delta$ is chosen small enough (depending on $C$).
	Thus by integration on $[0, t]$, for any $t \in\left[0, T_*\right]$, it holds 
 \[ D(t) \geqslant D(0)+\frac{\gamma g_0 t}{4 \alpha}. \]
	Therefore ,
	\begin{equation*}
	\frac 1{q({r}(t))} \geq \left(\frac {1- {C} \delta^2}{q({r}(0))}+ \frac{\gamma g_0 t}{4\alpha }\right)\big(1-{C}\delta^{2}\big) \ge \frac{1}{2} \delta^{-2}+\frac{\gamma g_0}{8\alpha}t ,
	\end{equation*}
	which in turn implies that
	\begin{equation*}
	q({r}(t))\le \frac{\frac{8\alpha}{\gamma g_0}}{t+\frac{4\alpha}{\gamma g_0}\delta^{-2}} \le \frac{1}{2}\frac{C_*}{t+\delta^{-2}}
	\end{equation*}
as soon as $C_* \ge 2$ is chosen such that $C_* \ge \frac{16\alpha}{\gamma g_0}$.
	This strictly improves the estimate of $q({r}(t))$ in \eqref{BS-2}.

 \bigskip
 
\bigskip
 
	\emph{Transversality condition.}

        Denote $\tilde b(t) = (t_*+\delta^{-1})^{2} b(t)$.
        Assume that there exists $t_* \in [0,\TS]$ such that $\tilde b(t_*)=1$.
        We claim that for this specific $t_*$, $\dot{\tilde b}(t_*) >0$.

        For this, we use \eqref{eq:b}, which gives
        \[
	\dot b(t) \ge 2\nu^{+}b(t)-C(\mathcal{N}^3+ q_*(r(t))\mathcal{N}).
	\]
        Due to the bootstrap bounds \eqref{BS-2} on $\q N(t) \le C_* \delta$ and $q(r(t)) \le C_*\delta^2$, there hold
        \begin{align*}
        \dot{\tilde b}(t_*) & = (t+\delta^{-1})^2 \dot b(t_*) +  2 (t_*+\delta^{-1}) b(t_*) \\
        & \ge 2 \nu^{+} \tilde b(t_*) + \frac{2 \tilde b(t_*)}{t+\delta^{-1}} -C(t_*+\delta^{-1})^2(\mathcal{N}^3(t_*)+q(r(t_*))\mathcal{N}(t_*)) \\
        & \ge 2 \nu^{+}  - 2C C_*^3 \delta \ge \nu^+ >0,
	\end{align*}
        provided $\delta >0$ is small enough (depending on $C_*$ and $C$ and $\nu_+>0$)

 \bigskip

 \emph{Step 2. A global solution via a topological argument}

We now consider all initial data $\bm {\gh a}^+ \in \q B_{\m R^K}(\delta^{5/4})$, and we dispose of a solution $\vec u = \vec u [\bm {\gh a}^+(0)]$ for which the bootstrap assumption \eqref{BS-2} holds on a maximal interval $[0,T_*(\bm {\gh a}^+(0))]$.

The above discussion performed for any given initial data of the form \eqref{def:ini} will allow us to justify the existence of at least one $K$-tuple 
	\[ \mathfrak{a}_*^{+}\in \mathcal{B}_{\mathbb{R}^{K}}(\delta^{5/4}) \]
 such that $\TS(\bm {\gh a}^+_*)=\infty$.
 
 We argue by contradiction and assume for its sake that, for all $\bm {\gh a}^+ \in \q B_{\m R^K}(\delta^{5/4})$, there hold $T_{*}(\bm a^+(0))<\infty$.

This assumption allows us to consider the map
\[ \Phi: \bar{\q B}_{\m R^K}(\delta^{5/4}) \to \m R^K, \quad \bm {\gh a}^+ \mapsto (T_*(\bm {\gh a}^+) +\delta^{-1}) \bm a^+(T_*(\bm {\gh a}^+)). \]

Let $\bm {\gh a}^+ \in  \bar{\q B}_{\m R^K}(\delta^{5/4})$. As $T_*(\bm {\gh a}^+) < +\infty$, due to the discussion in \emph{Step 1}, we have bootstrap saturation for the condition on $b$. This shows that $|\Phi(\bm {\gh a}^+)| = \sqrt{\tilde b(T_*(\bm {\gh a}^+))} = 1$. Hence $\Phi: \bar{\q B}_{\m R^K}(\delta^{5/4}) \to \m S_{\m R^K}$ takes values in the unit sphere.

Furthermore, the transversality condition shows that $\bm{\gh a}^+ \mapsto T_*(\bm {\gh a}^+)$ is continuous. By continuity of the flow of \eqref{DLKG} (and of the modulation), we infer that $\Phi$ is continuous.

Invoking again the transversality condition, we see that there is an instantaneous exit on the boundary: $T_*(\bm{\mathfrak a}^+) =0$ for all $\bm {\gh a}^+ \in \m S_{\m R^K}(\delta^{5/4})$. Therefore, $\Phi$ is a retraction in the sense that
\[ \Phi|_{\m S_{\m R^K}(\delta^{5/4})} = \delta^{-5/4} \Id. \]
The existence of such a map $\mathcal{M}$ contradicts Brouwer's no retraction theorem for continuous maps from the ball to the sphere.

We conclude that there exists at least one  $\bm{\mathfrak{a}}_*^{+}\in \mathcal{B}_{\mathbb{R}^{K}}(\delta)$ such that $T_{*}(\bm{\mathfrak a}_*^+)=\infty$.

\bigskip

\emph{Step 3. Precise description of the asymptotics}

 For the specific $\bm{\mathfrak{a}}_*^{+}$ obtained in \emph{Step 2}, denote by $u_\sharp$ the associated solution given by \eqref{def:ini}. It is globally defined for positive times;  the bootstrap bounds \eqref{BS-2} hold for all times $t \ge 0$, and the bounds \eqref{eq:drk2} and \eqref{est:rot} as well. From there, equations \eqref{est:translation} and  \eqref{eq:drk2} can be rewritten
\begin{gather*}
\dot \tau_\sharp(t) = O((t+ \delta^{-1})^{-\theta})
\dot r(t) = \frac{\gamma g_0}{2\alpha} q(r(t)) + O((t+ \delta^{-1})^{-\theta}).
\end{gather*}
By integration, we infer that $\tau_\sharp$ admits a limit $\tau_\infty$ as $t \to +\infty$, and
\[ \tau_\sharp(t) - \tau_\infty = O((t+ \delta^{-1})^{1-\theta}). \]
We can also integrate \eqref{equ:dis-sym} (see also \cite{CMYZ}), as follows. We first claim the convergence
\begin{equation} \label{conv:q}
\frac{1}{t q(r(t))} \to \frac{\gamma g_0}{2\alpha} \quad \text{as} \quad t \to +\infty.
\end{equation}
For this, let us first prove that for some constant $\tilde C = \tilde{C}(\delta)$,
\begin{equation}\label{limit exist}
\left|\frac{1}{q(r(t))}-\frac{\gamma g_0}{2 \alpha} t\right| \lesssim \frac{t}{|\ln \delta|}+\tilde{C}(\delta)
\end{equation}

Indeed, define $\ds E(t)=\frac{1}{q(r(t))}\cdot e^{C \q M (t)}$, and using \eqref{equ:dis-sym} and \eqref{eq:Mbis}, we infer
\begin{align}
e^{-C\mathcal{M}(t)} \dot E(t) &=\frac{d}{dt} \left[ \frac{1}{q(r(t))}\right] +\frac{{C}}{q({r}(t))} \dot{\mathcal M} \nonumber  \\
& \le \frac{\gamma g_0}{2\alpha}+C\frac{\mathcal{N}^2}{q(r(t))}+C\frac{\mathcal{N}^2}{q(r(t))}+C {r(t)}^{-1}-{C}q^2({r(t)}) \nonumber \\
& \le \frac{\gamma g_0}{2\alpha}+ \frac{2C}{r(t)}. \label{limit exist-2}
\end{align}
As $r(t) \gtrsim |\ln \delta|$ and  $|e^{C\q M(t)} -1| = O(C\delta^2)$, we infer from there that
\[
\frac{\mathrm{d}}{\mathrm{d} t} E(t)\leqslant \frac{\gamma g_0}{2\alpha}\left(1+\frac{O(1)}{|\ln \delta|}\right)\left(1+O\left(C \delta^2\right)\right) \leqslant \frac{g_0}{\alpha}\left(1+\frac{O(1)}{|\ln \delta|}\right),
\]
which is \eqref{limit exist}. Now, as $\q N(t), q(r(t)) \to 0$, we can repeat the computations in the bootstrap of \emph{Step 1} at later times, so that \eqref{limit exist} holds for any $\delta >0$, provided $t$ is large enough (depending on $\delta$). After dividing by $t$, this exactly boils down to \eqref{conv:q}. Using the expansion  \eqref{Qdec} of $q$, this is equivalent to
\[
\lim _{t \rightarrow \infty} \frac{1}{t} r(t)^{\frac{d-1}{2}} e^{r(t)}=\frac{\kappa\gamma g_0}{2\alpha} .
\]
Setting $c_0=\ln\left( \frac{\kappa \gamma g_0}{2\alpha}\right) $, this implies
\begin{equation}\label{equ:asyr}
r(t)=\ln t-\frac{1}{2}(d-1) \ln \ln t+c_0+s(t), \quad \text {where}  \quad s(t) \to 0 \text { as } t \to \infty .
\end{equation}
Let us now refine the convergence of $s$. For this, we go back to \eqref{eq:drk2}, which writes
\begin{align*}
\dot{r} & =\frac{\gamma}{2\alpha}g(r)+O\left(t^{-\theta}\right)=r^{-\frac{d-1}{2}} e^{-r+c_0} (1+o(1))+O\left(t^{-\theta}\right) \\
& =r^{-\frac{d-1}{2}} e^{-r+c_0}\left(1+O\left(\frac{1}{\ln t}\right)\right).
\end{align*}
Therefore,
\[ 
\frac{\mathrm{d}}{\mathrm{d} t}\left(r^{\frac{d-1}{2}} e^{r-c_0}\right)=1+O\left(\frac{1}{\ln t}\right) .
\]

By integration,
\begin{equation}\label{equ:dis-as}
r(t)-c_0+\frac{d-1}{2} \ln r(t)=\ln t+O\left(\frac{1}{\ln\ln t}\right) .
\end{equation}

Then from \eqref{equ:asyr} and \eqref{equ:dis-as},
we infer
\[
s(t)=-\frac{d-1}{2}\ln\left(\frac{r(t)}{\ln t} \right) +O\left( \frac{1}{\ln t}\right) =\left( \frac{d-1}{2}\right)^2\frac{\ln \ln t}{\ln t}+O\left( \frac{1}{\ln t}\right) 
\]
This shows that
\[ r(t) = \ln t-\frac{1}{2}(d-1) \ln \ln t+c_0 + O\left( \frac{\ln \ln(t)}{\ln t}\right). \]
Finally, arguing as for \eqref{est:rot}, we see that for all $k=1, \dots, K$,
\[ \dot R_\sharp(t) \omega_k = O((t+\delta^{-1})^{-\theta}). \]
Therefore, $R_\sharp|_{\spn \Omega}$ admits a limit $R_\infty \in \mc L(\spn(\Omega), \m R^d)$ which is an isometry, which we can extend into an element of $O_d(\m R^d)$ still denoted $R_\infty$. Then
\[ \forall k=1, \dots, K,  \quad R_\sharp(t) \omega_k - R_\infty \omega_k = O((t+\delta^{-1})^{1-\theta}). \]
Hence (allowing constants to depend on $\delta$):
\begin{align*}
z_k(t) - \tau_\infty & = y_k(t)   + \frac{\ell_k(t)}{2\alpha} - \tau_\infty=  \lambda_\Omega r(t) R_\sharp(t) \omega_k +  \frac{\ell_k(t)}{2\alpha} \\
& = \lambda_\Omega r(t) R_\infty \omega_k + O(t^{1-\theta} \ln(t))
\end{align*}
and in conclusion,
\[ \vec u_\sharp = S[\bm z, 0, \bm \sigma] + \vec \e(t) \quad \text{with} \quad \| \vec \e(t) \|_{H^1 \times L^2} = O(t^{-1}). \qedhere\]
\end{proof}

\begin{proof}[Proof of Theorem \ref{thm:abs}]
Apply Proposition \ref{pr:exist} to the rescales $C \Omega$ for some sufficiently large $C$, so that $q_*(C \Omega) \le \delta$, with $\bm \ell =0$ and $\vec \e_\perp =0$: it provides a solution $u_\sharp$ (together with $\tau_\infty$, $R_\infty$) defined globally for positive times and satisfying \eqref{conv:prop} and \eqref{conv:prop2}. Now consider 
\[ u(t,x) : = u_\sharp(t, R_\infty^{-1}(x-\tau_\infty)). \]
This is a solution of \eqref{DLKG} which satisfies
\[ \vec u(t) = \vec S[\bm z, 0 ,\bm \sigma] + \vec \e(t), \]
with $\| \vec \e(t) \|_{H^1 \times L^2} \lesssim t^{-1}$ and (as $\lambda_{C\Omega} = \lambda_{\Omega}/C$), for $k=1, \dots, K$,
\[ z_k(t)  = \lambda_{C\Omega} r(t) C \omega_k = \lambda_\Omega  \left( \ln (t)  - \frac{d-1}{2} \ln \ln t + c_0 \right) \omega_k + O \left( \frac{\ln \ln t}{\ln t} \right). \]
Therefore, $u$ is the desired solution, with $c_\Omega = \lambda_\Omega c_0 = \lambda_\Omega \ln \left( \frac{\kappa \gamma g_0}{2\alpha} \right)$. 
\end{proof}

The proof of Corollary \ref{cor:sym} follows immediately, as in the cases mentioned, the rigidity of $(\Omega, G)$ is clear. It remains to consider the orthoplex without a center, which we do now.

\begin{proof}[Proof of Corollary \ref{cor:orthoplex}]
The facts that $G$ is subgroup of the symmetry group of $\Omega$ and that $\sigma$ is an homomorphism on $G$ are clear.

Therefore, it suffices to see that $(\Omega,G)$ is a rigid configuration. Consider 
\[G_0 = \{ S^k R_1^{\alpha_1} \cdots R_d^{\alpha_d} : k \in \llbracket 0, d \rrbracket, \alpha_j \in \{ 0,1 \} \}, \]
where $R_j$ is the reflection across ${\bf e}_j^\perp$: $R_j {\bf e_i} = {bf e_i}$ for $i \ne j$ and $R_j {\bf e_j} = - \bf e_j$, and $S$ is the cyclic permutation $S {\bf e_j} = {\bf e} _{j+1}$, with the convention that ${\bf e_{d+1}} = {\bf e_1}$. Observe that $G_0$ is a subset of $G$, and that $(\Omega,G_0)$ is already rigid. Indeed, if $\Omega'$ is a subset of $\m R^d$ of cardinal $2d$ invariant by $G_0$, consider $\bs x = (x_1, \dots, x_d) \in \Omega'$. Up to applying reflection, we can assume that all $x_j \ge 0$, and up to applying a rotation (which would be the one appearing in the definition of rigidity), we can assume that $x_1 \ge x_2 \ge \cdots \ge x_d$. As $2d < 2^d$, $x_d=0$ (otherwise all the $(\pm x_1, \dots, \pm x_d) \in \Omega'$, and they are too numerous). Assume $\bs x \ne 0$, that is $x_1 >0$. Letting $S$ act we see that all $S^k \bs x$ are different for $0 \le k \le d-1$ (except if $\bs x=0$), and if $x_2 >0$, we obtain at least $4d$ points in $\Omega'$, which too much: indeed, the reflections give that the 4 distinct points $(\pm x_1, \pm x_2, x_3, \dots, x_d)$ belong to $\Omega'$, and the action of the cyclic permutation give $d$ points for each of them, as one can be convinced by considering the last non zero coordinate.

Hence $\bs x = \lambda \bs e_1$, $\Omega' = \lambda \Omega$ and $\phi = \lambda \Id$ fits. If $0 \in \Omega'$, the same argument shows that is has cardinality $1$ or at least $2d+1$, both of which are not allowed. Hence $(\Omega, G_0)$ is rigid, and a fortiori, $(\Omega, G)$ as well.
\end{proof}

\section{Non existence of same sign multi-solitons} \label{sec:non}

In this section, we prove Theorem \ref{thm:non-existence}. We start with a (conditional) stability result.

\begin{proposition}[Long time convergence] \label{non existence}
Let $K \ge 2$ and $\vec{u}$ be a $K$-soliton of \eqref{DLKG} who\-se solitons all have the same sign.

Let $\delta >0$ small enough, and $T_\delta \in \mathbb{R}$ such that $\vec{u}(t)$ admits a decomposition as in Lemma \ref{le-dec} in a neighbourhood of $T_\delta$ with
\begin{equation}\label{starting point}
\mathcal{N}\left(T_\delta\right) \leqslant \delta ,\quad q_*\left(\bm z(T_\delta) \right) \leqslant \delta^2,\quad b(T_\delta)\le \delta^2.
\end{equation}
Then this decomposition holds  on $[T_\delta, \infty)$ and 
\begin{equation} \label{non-exis:decay}
 \forall t \ge T_\delta, \quad q_*(\bm y(t)) \lesssim \mathcal{N}(t)^2  \lesssim \delta^2.
\end{equation}
\end{proposition}

\begin{proof}
For $\delta >0$ small, the existence of $T_\delta$ satisfying \eqref{starting point} is a consequence of definition \ref{def-multi wave} and lemma \ref{le-dec} (and continuity of the flow).

Let $C_* >0$ be a large constant to be fixed later. We introduce the following bootstrap estimates
\begin{equation}\label{BS-non}
\mathcal N(t) \le C_*\delta, \qquad b(t)  \le C_* \delta^2,
\end{equation}
and we set
\begin{equation}\label{bs-non}
\TS= \sup\left\{ t\in [T_\delta,\infty)\text{ such that \eqref{BS-non} holds on } [T_\delta,t] \right\}.
\end{equation}
By continuity of the flow, $\TS > T_\delta$. We argue by contradiction, and assume that $T_* <+\infty$. We now work on $t \in [T_\delta, T_*)$.

We consider Lemma \ref{ene expand}. As $\vec u$ is a $K$-soliton, $E(\vec u(t)) \to K E(Q)$ as $t \to +\infty$, so that the energy identity writes
\[ 2 \alpha \int_t^{+\infty} \| \partial_t u (t) \|_{L^2}^2 \ud t = E(\vec u(t)) - K E(Q). \]
Notice that in \eqref{exp:energy}, all $\sigma_j \sigma_k =1$ as the sign are all the same. Therefore, \eqref{exp:energy} gives:
\[ \int_{t}^\infty \| \partial_t u\|_{L_x^2}^2 \ud t + c_1 g_0 \sum_{j \ne k} q(|z_j - z_k|) = O \left( \frac{q_*(\bm z(t))}{|\ln(q_*(\bm z(t))|} + \q N(t)^2 \right). \]
As $\sum_{j \ne k} q(|z_j - z_k|) \ge q_*(\bm z(t))$, we get
\[ \int_{t}^\infty \| \partial_t u\|_{L_x^2}^2 \ud t + q_*(\bm z(t)) \lesssim  \frac{q_*(\bm z(t))}{|\ln(q_*(\bm z(t))|} + \q N(t)^2. \]
Choosing $\delta>0$ small enough allows us to absorb the first term of the right hand side, and as $q_*(\bm z(t)) \sim q_*(\bm y(t))$, we arrive at
\begin{equation}\label{ene-expand-+}
\forall t \in [T_\delta, T_*), \quad \int_{t}^\infty \| \partial_t u\|_{L_x^2}^2 \ud t + q_*(\bm y(t)) \lesssim \q N(t)^2.
\end{equation}
In particular, estimates \eqref{eq:b}, \eqref{eq:damped1} and \eqref{eq:damped2} now write on $[T_\delta, T_*)$
\begin{align} 
|\dot b - 2\nu_+ b| & \lesssim \q N^3, \label{eq:b_2} \\
\frac{\ud}{\ud t} \q F + 2\mu \q F & \lesssim \q N^3, \label{eq:F_2} \\
\frac{\ud}{\ud t} \q G + 2\mu \q G & \lesssim \q N^3. \label{eq:G_2}
\end{align}

We now proceed to improve the bootstrap estimate \eqref{BS-non}. The inspiration is similar to the proof Proposition \ref{pr:exist}, but the computations are now done in a context of decay estimates. Again, $C_* \ge 1$ will be chosen large and then $\delta >0$ small with in particular the bound $C_*^3 \delta \le 1$, and in the following computations, implicit constants do not depend on $\delta$ nor $C_*$.

\bigskip

\emph{Estimate on $\q N$.}

From~\eqref{eq:damped1} and, it holds on $[\TD,\TS)$,
\begin{align*}
	\frac{d}{dt}\left[ e^{2\mu t}\mathcal F\right] &=e^{2\mu t}\left(\frac{d}{dt}\mathcal F+2\mu \mathcal F \right) \lesssim e^{2 \mu t} \mathcal N^3 \\
	&\lesssim e^{2\mu t} C_*^3  \delta^3 \lesssim e^{2\mu t}\delta^{2} ,
	\end{align*}
	for $\delta$ small enough. Thus, integrating the above estimate on $[0,t]$,
	for any $t\in [0,\TS)$, 
	\begin{equation*}
	\mathcal F(t)-e^{-2\mu (t-\TD)}\mathcal F(\TD) \lesssim \delta^{2}\int_{\TD}^{t} e^{-2\mu(t-s)} \ud s \lesssim \delta^{2}.
 \end{equation*}
	From~\eqref{eq:new2} and \eqref{BS-2}, we conclude that
	\[ \mathcal F(t)\lesssim \delta^{2}. \]
	In view of \eqref{eq:new2}, we obtain
	\begin{equation*}\label{qF}
	\|\vve\|_\ENE^2 \lesssim \mathcal F+b
	\lesssim C_* \delta^{2}
	\end{equation*}
	Arguing similarly for the quantity $\BB$, we infer
	\begin{equation} \label{G-new}
	 |\bm \ell|^2 +  | \bm a^-|^2 \lesssim \BB \lesssim \delta^{2}.
	\end{equation}
	Hence we obtain, for all $t \in [\TD, \TS)$,
	\begin{equation}\label{N-new}
	\mathcal N(t) \lesssim \sqrt{C_*}\delta.
	\end{equation}
	Fix $C_*$ large enough to absorb the implicit constant, and then for $\delta>0$ small, this strictly improves the bootstrap estimate \eqref{BS-2} on~$\mathcal N$ on the interval $[0,\TS)$.

	\bigskip
	
	\textit{Estimate on $b$.}
	
	As $\partial_t u=\eta-\sum_{k=1}^K \left(\ell_k \cdot \nabla\right) Q_k$ we have 
	\[ \forall t \in [\TD,\TS), \quad \|\eta(t) \|_{L^2}^2 \lesssim\left\|\partial_t u(t) \right\|_{L^2}^2+|\bm \ell| ^2(t). \] 
	By the definition of $a_k ^{ \pm}$, one has
\[
	a_k^{+}=\zeta^{+}\left\langle \varepsilon, Y_k \right\rangle+\left\langle \eta, Y_k\right\rangle, \quad a_k^{-}=\zeta^{-}\left\langle \varepsilon, Y_k\right\rangle+\left\langle\eta, Y_k \right\rangle,
\]
	and thus
\[
	a_k^{+}=\frac{\zeta^{+}}{\zeta^{-}} a_k^{-}+\frac{\zeta^{-}-\zeta^{+}}{\zeta^{-}}\left\langle \eta, Y_k \right\rangle .
\]
Therefore,
\begin{equation}\label{b-improve}
	b(t) \lesssim \|\eta\|_{L^2}^2+|\bm a^-|^2 
	\lesssim
	\|\partial_t u\|_{L^2}^2+|\bm \ell| ^2+|\bm a^-|^2 
	\lesssim 
	\|\partial_tu \|_{L^2}^2+\q G.
\end{equation} 

As $T_* <+\infty$, and in view of the improved estimates of $\q N$, $b(t) \to C_* \delta^2$ as $t \to T_*$. Hence, by continuity of $b$, there exists $t_2 \in\left[T_\delta, T_*\right)$ such that
\[ 
	b\left(t_2\right)=\frac{C_*}{2} \delta^2, \quad \text{and for all } t \in \left[T_\delta, t_2\right), \quad b(t)<\frac{C_*}{2} \delta^2. \]	
Still by continuity, there exist $t_1 \in \left[T_\delta, t_2\right)$ such that
\[ 
b\left(t_1\right)=\frac{C_*}{4} \delta^2 \quad \text { and for all }  t \in [t_1,t_2], \quad b(t)>\frac{C_*}{4} \delta^2.
\]
Now compute, using	\eqref{eq:b_2}
	\begin{align*}
	\dot b(t) & = 2 \nu_+ b + O(\q N^3) = 2 \nu_+ b +O(C_*^3\delta^3).
	\end{align*}
 and so, on $[t_1,t_2]$,
\[ \left| \frac{\dot b}{2\nu_+ b} -1 \right| \lesssim C_*^2 \delta \le \delta^{1/2}, \]
if $\delta>0$ small enough. Integrating over $[t_1,t_2]$ we get that
\[
t_2-t_1 = \frac{\ln 2}{2 \nu^{+}}+O\left( \delta^{1/2} \right).
\]
Then
	\begin{equation}\label{bge}
	\int_{t_1}^{t_2} {b}(t) \ud t \gtrsim C_* \delta^2.
	\end{equation}
	However, due to \eqref{b-improve} and the improved bound \eqref{N-new} on $\q N$ (and on $\q G$) ,
	\begin{align*}
	\int_{t_1}^{t_2} b(t) \ud t &\lesssim \int_{t_1}^{t_2} \| \partial_t u\|_{L^2}^2\ud t +\int_{t_1}^{t_2} \q G(t)\ud t \\
	& \lesssim \int_{T_{\delta}}^{+\infty} \| \partial_t u\|_{L^2}^2\ud t + \delta^2 \lesssim \q N(T_\delta) + \delta^2 \lesssim \delta^2.
	\end{align*}
Choosing $C_*$ large enough (depending on the implicit constants), this gives a contradiction with \eqref{bge}.

Therefore $T_* =+\infty$, and the bounds \eqref{BS-non} and \eqref{ene-expand-+} hold for all times $t \ge T_\delta$. 


\end{proof}

We are now in a position to conclude the proof of Theorem \ref{thm:non-existence}.
\begin{proof}[Proof of Theorem \ref{thm:non-existence}]
	
	We argue by contradiction and assume that $\vec u$ is a $K$-soliton where all the solitons have the same sign ($K \ge 2)$. We can apply the previous Proposition \ref{non existence}: there exist $T_\delta$ such that for $t \ge T_\delta$, $\vec u$ admits a decomposition in the sense of Lemma \ref{le-dec} and satisfies the bound \eqref{non-exis:decay}.

 \bigskip
	
\emph{Step 1: Lower bound on $q_*$}
	
	Recall the ODE of the centers \eqref{eq:dist}:
	\[ \dot y_k = \frac{1}{2\alpha} \sum_{j \in \llbracket 1, K \rrbracket \setminus \{ k \}}  \frac{y_{j}-y_k}{|y_{j}-y_k|} g(|y_{j}-y_k|)+O(\mathcal N^2+q_*(\bm y)^\theta) = O(\q N^2 + q_*(\bm y)). \]
	Hence for $j \ne k$, denoting $r_{j,k} = |y_j - y_{k}|$, we infer
	\begin{align*}
	\dot r_{j,k} & = \frac{(\dot y_j - \dot y_{k}) \cdot ( y_j - y_{k})}{r_{j,k}} = O(\mathcal N^2+q_*(\bm y)).
	\end{align*}
	And so,
	\begin{align}
	\frac{\ud }{\ud t} \frac{1}{q_*(\bm y)} &= -\frac{1}{q_*^2(\bm y)}  \sum_{j \ne k} q'(r_{j,k}) \dot r_{j,k}  \nonumber \\
	& = \frac{1}{q_*^2(\bm y)}  \sum_{j \ne k} O \left(  \frac{q(r_{j,k})}{r_{j,k}} (\q N^2 + q_*(\bm y)) \right) \nonumber  \\
	& \le C \left( \frac{\q N^2}{q_*(\bm y)} +   \sum_{j \ne k} \frac{1}{r_{j,k}}   \right). \label{der:q_*}
	\end{align}
	
\emph{Step 2: Estimate on the refined distance $R_*$}
	
	If we define $\ds R_*(t)=\frac{1}{q_*(\bm y)} e^{C\q M}$,  then by \eqref{eq:Mbis} and \eqref{der:q_*}
	then we get 
	\begin{align*}
	\frac{\ud}{\ud t}R_*(t) &=e^{C\q M}\left( \frac{\ud}{\ud t}\frac{1}{q_*(\bm y)}+C\frac{1}{q_*}\frac{\ud}{\ud t} \q M\right) 	 \\
	& \le C \frac{\mathcal{N}^2}{q_*(\bm y)}-C\frac{\mathcal{N}^2}{q_*(\bm y)}+C\sum_{j \ne k}r^{-1}_{j,k} \le \frac{C'}{|\ln \delta|}.
	\end{align*}
	By integration, we infer that for $t \ge T_\delta$,
	\[ R(t) \le R(T_\delta) + \frac{C}{|\ln \delta|} (t-T_\delta). \]
	Now $R(T_\delta) \lesssim \delta^{-2}$ and $\q M = (\delta^2)$, so that
	\[ \frac{1}{q_*(\bm y(t))} \lesssim \delta^{-2} + \frac{1}{|\ln \delta |}(t-T_\delta), \]
	and so
	\[
	q_*(\bm y(t))\gtrsim \frac{1}{\delta^{-2}+|\ln \delta|^{-1}(t-T_\delta)}
	\]

\bigskip

\emph{Step 3: Contradiction}	
	
Due to \eqref{non-exis:decay}, we infer from the above bound that
\begin{equation} \label{est:N_div}
\int_{T_\delta}^{\infty} \q N(t)^2 \ud t \gtrsim \int_{T_\delta}^{\infty} q_*(\bm y(t)) \ud t = +\infty.
\end{equation}
Now recall \eqref{eq:Mbis}: in view again of \eqref{non-exis:decay}, $q_*(\bm y)^2 \lesssim \q N^4$ and so, for $t \ge T_\delta$
\[ \frac{\ud \q M}{\ud t} \le - \frac{1}{2} \q N^2. \]
Integrating in time, we infer that $\q M(t) \to -\infty$ as $t \to +\infty$. But this contradicts that
\[ |\q M(t)| \lesssim \q N^2 \lesssim \delta^2 \]
is bounded (and small). This concludes the proof.
\end{proof}

\section{Construction on the line in dimension \texorpdfstring{$d\ge2$}{d >= 2}} \label{sec:line}

The goal of this section is to prove Theorem \ref{th:2}, that there exists a $K$ soliton whose centers lie on a line. For this, as mentioned in the introduction, we work with data with cylindrical symmetry: by performing a rotation, we can assume that symmetry axis is $e =e_1 = (1,0, \dots, 0) \in \m R^d$, and so we ask that the dependence on the variables $x_2, \dots x_d$ is only through $\rho = \sqrt{x_1^2 + \cdots + x_d^2}$.

This corresponds to the symmetry group made of matrices
\[  \begin{pmatrix}
    1 & 0 & \cdots & 0 \\
    0 \\
    \vdots & & M \\
    0
    \end{pmatrix} \quad \text{with} \quad M \in O_{d-1}(\m R).
\]
This symmetry is preserved by the flow. We denote $\mathcal {H}_\mathcal{C} $ the space of data with cylindrical symmetry. 


We argue as in Section 3 and prove an analogue of Proposition \ref{pr:exist}, as follows.

\begin{proposition}\label{prop:1d constr}
For $\delta > 0$ small enough let any
\[ \left\{
\begin{array}{l}
(\ell_{k}(0))_{k=1, \dots, K} \in \mathcal{B}_{\mathbb{R}^{K}}(\delta),
(z_{k}(0))_{k=1, \dots, K} \in \mathbb{R}^{K} \\
\text { satisfying } \min \left\{ z_{k+1}(0) -z_{k}(0)  , k =1, \dots, K-1 \right\}>5|\log \delta|, \\
\vec{\varepsilon}_{\perp}(0) \in \mathcal{B}_{H^1 \times L^2}(\delta) \cap \mathcal{H}_{\mathcal{C}},\quad
 \left\langle\vec{\varepsilon}_{\perp}(0), \vec{Z}_{k}^{+}(0)\right\rangle=0 \text { for } k =1, \dots, K
\end{array}\right.\]
There exists $\bm{\gh a}^{+}(0)=\left(\gh a_{k}^{+}(0)\right)_{k =1,\dots, K } \in \bar{\mathcal{B}}_{\mathbb{R}^{K}}\left(\delta^{\frac{5}{4}}\right)$ such that the solution $\vec{u}_{\sharp}$ of \eqref{DLKG} with the initial data
\[ \vec{u} (0)= S[\bm z(0), \bm \ell(0), \bm \sigma] +\vec{W}\left(\bm{\gh a}^{+}(0)\right)+\vec{\varepsilon}_{\perp}(0) \]
is global, with cylindric symmetry and admits a decomposition
\[ \vec u_\sharp(t) = S[\bm z(t), 0, \bm \sigma] + \vec \e(t) \]
where, for $k =1, \dots, K$ 
\begin{equation} \label{eq:zk}
z_k(t)=\left(k-\frac{K+1}{2}\right)\left(\ln t-\frac{d-1}{2}
\ln\ln t \right) + O \left( 1 \right),
\end{equation}
and $\| \vec \e(t) \|_{H^1 \times L^2} = O(t^{-1})$.
\end{proposition}

Notice that as $G$ is no longer rigid, because all points on $\m R e_1$ are invariant, one should allow modulation along $\m R e_1$.

\bigskip

At leading order, the soliton centers solve the ODE system for the unknowns $\zeta_1, \dots, \zeta_K$ (with values in $\m R$):
\begin{equation} \label{def:ODE_sys}
\dot \zeta_k = \frac{1}{2\alpha} (g(\zeta_k -\zeta_{k-1}) - g(\zeta_{k+1}- \zeta_k))
\end{equation}
(with the conventions that 
\begin{equation} \label{eq:conv}
\zeta_0 = -\infty,  \quad \zeta_{K+1} = +\infty \quad \text{and} \quad g(+\infty) = 0.
\end{equation}
Before we proceed with Proposition \ref{prop:1d constr}, let us first state a stability result concerning forward global solutions of the ODE system \eqref{def:ODE_sys}, with a perturbation.

\begin{proposition} \label{prop:stab}
Let $y_1, \dots, y_K \in \mc C^1([0,+\infty),\m R)$ satisfy for all $k =1, \dots,K$ and as $t \to +\infty$
\[ \dot y_k =  \frac{1}{2\alpha} (g(y_k -y_{k-1}) - g(y_{k+1}- y_k)) + O \left( \sum_{j=1}^{K-1} g(y_{j+1} -y_{j})^{\theta} \right). \]
(with the convention \eqref{eq:conv}) and $y_{k+1}(t) - y_k(t) \to +\infty$ for all $k=1, \dots, K-1$.

Then for all $k =1, \dots, K$, and as $t \to +\infty$,
\begin{equation} \label{expansion_yk}
y_k(t) = \left( k - \frac{K+1}{2} \right) \left( \ln t - \frac{d-1}{2} \ln \ln t \right) + O(1).
\end{equation}
\end{proposition}

\begin{remark} \label{nb:stab_ODE}
One could in principle improve the $O(1)$ in \eqref{expansion_yk} and in \eqref{eq:th:z} as well, to a control of the same type as in \eqref{dis:new}, namely $\tau_k + O\left( \frac{\ln \ln t}{\ln t} \right)$, for some universal constant $\tau_k$. This could be done by revisiting the asymptotic stability arguments developed in \cite{CZ13}.
As the proof of Proposition \ref{prop:stab} is already technical, we will not pursue this issue here.
\end{remark}

\bigskip

We postpone the proof of Proposition \ref{prop:stab} to the appendix, and using it, we complete that of Proposition \ref{prop:1d constr}.

\begin{proof}[Proof of Proposition \ref{prop:1d constr}]
Given $\bm {\gh a}^+ = (a_\omega^+)_{1 \le k \le K} \in \q B_{\m R^{K}}(\delta^{5/4})$, we consider the solution $\vec u(t)$ of \eqref{DLKG} with initial data
\begin{align} \label{def:ini2}
\vec u(0) = \vec S[\bm z(0), \bm \ell(0), \bm \sigma] + \vec W(\bm {\gh a}^+) + \vvep(0).
\end{align}

\emph{Decomposition.}
For any $T > 0$ such that the solution $\vec u$ satisfies on $[0,T]$,
\[ \vec u(0) =  S[\bm \zeta(t), 0, \bm \sigma] \]
for some continuous function $\bm \zeta(t)$, we decompose it according to Lemma~\ref{le-dec} taking into account cylindrical symmetry: for $t \in [0,T]$ we therefore dispose of $\bm z(t), \bm \ell(t) \in \m R^K$, $\vec \e(t) \in \q H_{\q C}$, and from there, $\bm y(t) \in \m R^K$ and $\bm a^\pm(t) \in \m R^K$. 

As simple continuity argument shows that the coordinates of $\bm z$  and $\bm y$ remain well separated and ordered, so that for all $k=1, \dots, K-1$ and $t \in [0,T]$,
\[ z_{k+1}(t) -z_k(t), \ y_{k+1} (t) - y_k(t)  \gtrsim \ln| \delta|, \]
so that for $k=0,\dots K$, we define
\[ r_k(t) = y_{k+1}(t) - y_k(t) \ge 0. \]
(The bounds above will be much improved below in the bootstrap argument, but it is a convenient remark to make now for notations).

As in Proposition \ref{pr:exist}, due to the properties of the function $W$ in Lemma~\ref{le-W} and the orthogonality properties \eqref{ortho} of  $\vvep(0)$ assumed in~\eqref{eq:pr:10}, the initial data $\vec u(0)$ is modulated, in the sense that $ (\bm{z}(0),\bm{\ell}(0))$ and
\begin{equation*}
    \vve(0)=\vec W(\bm {\gh a}^+) + \vvep(0),
\end{equation*}
are the parameters of the decomposition of $\vec u(0)$. Lemma~\ref{le-W} applies and, for $k=1, \dots, K$,
	\begin{equation*}
a_k^+(0) =  \langle \vve(0),\vec Z_{k}^+(0)\rangle =  \langle \vec W(\bm {\gh a}^+), Z_k^+ \rangle =  \gh a_k^+.
\end{equation*}
We slightly adapt the function estimating the distance and define
\[ q_*(\bm y(t)) = \sum_{k=1}^{K-1} q(y_{k+1}(t) - y_k(t)). \]
We now define $\q F, \q G, \q N, \q M$ as in Section \ref{S:2.4}. Then at initial time $t=0$, it holds from~\eqref{eq:pr:10}
\begin{equation}\label{at:02}
    \mathcal N(0)\lesssim \delta,\quad q_*(\bm z(0))\lesssim \delta^2 \quad \text{and} \quad b(0) \le \delta^{5/2}.
\end{equation}
Also, the estimates \eqref{eq:new2}, \eqref{eq:b}, \eqref{eq:damped1} and \eqref{eq:damped2}, \eqref{eq:Mbis} hold, let us recall them here:
\begin{align}
|\dot b - 2 \nu^+ b| & \lesssim \q N^3 + q_*(\bm y), \label{eq:dot_b}\\
\frac{d}{dt} \q F + 2\mu \q F & \lesssim \q N^3 + \q N q_*(\bm y), \\
\frac{d}{dt} \q G + 2\mu \q G & \lesssim \q N^3 + \q N q_*(\bm y), \\
\frac{d}{dt} \q M + \q N^2 & \lesssim q_*(\bm y)^2. \label{eq:M_lyap}
\end{align}
The system \eqref{eq:dist}  for $\bm y$ now takes the form: for $k=1,\dots, K$,
\begin{equation} \label{eq:ode_sys} 
\dot y_k = \frac{1}{2\alpha} \left( g(y_{k} - y_{k-1}) - g(y_{k+1} - y_k) \right) + O(\q N^2 + q_*(\bm y)^\theta), 
\end{equation}

(with the convention \eqref{eq:conv}).

\bigskip

We also need a control on the evolution of $q_*(\bm y)$, as follows

\begin{claim} \label{cl:dot_q*}
There exists $\lambda >0$ such that for all $t \in [0,T_*]$
 \begin{equation}\label{le:q_* control} 
  \lambda - \frac{1}{\lambda q_*}\left(\mathcal{N}^2 + \sum_{k=1}^{K-1} r_k^{-1}q(r_k)\right) \le \frac{\mathrm{d}}{\mathrm{d}t}\left[\frac{1}{q_*}\right] \le \frac{1}{\lambda} + \frac{1}{\lambda q_*}\left(\mathcal{N}^2 + \sum_{k=1}^{K-1} r_k^{-1}q(r_k)\right).
\end{equation}
\end{claim}

\begin{proof}
Recall that $g(r)=g_0 q(r)+O(r^{-1}q(r))$. Due to \eqref{eq:ode_sys},\begin{align*}
\dot{r}_{k} & =\dot{y}_{k+1}-\dot{y}_{k} =\frac{1}{2 \alpha}( - g(r_{k+1})  +2 g(r_k) - g(r_{k-1})) + O( \q N^2 + q_*(\bm y)^\theta) 
\\
&= \frac{g_0 }{\alpha } ( - q(r_{k+1})  +2 q(r_k) - q(r_{k-1})) +O\left(\mathcal{N}^2+\sum_{k=1}^{K-1}r_k^{-1}q(r_k)\right).
\end{align*}
Therefore, recalling that $q'(r)=-q(r)+O(r^{-1}q(r))$, we obtain
\begin{align} \label{eq:est_dot_q*}
\dot q_*(\bm y) = \sum_{k=1}^{K-1}\dot{r}_kq'(r_k) &=-\frac{ g_0}{2\alpha}S+O\left(q_*(\bm y) (\mathcal{N}^2+\sum_{k=1}^{K-1}r_k^{-1}q(r_k))\right),
\end{align}
where 
\[ S=\sum_{k=1}^{K-1} \left( -q(r_{k})q(r_{k-1}) + 2q^2(r_k) -q(r_k)q(r_{k+1}) \right) =  q(\bm r)^\top A  \ q( \bm r) \]
and we denoted
\[ q(\bm r) = \begin{pmatrix}
q(r_0) \\
\vdots \\
q(r_K)
\end{pmatrix} \in \m R^{K+1} \quad \text{and} \quad A = \begin{pmatrix}
2 & -1 & 0 & \ldots & 0 \\
-1 & 2 & -1 & \ddots & \vdots \\
0 &-1 & 2 & \ddots & 0 \\
\vdots & \ddots & \ddots & \ddots & -1 \\
0 & \ldots& 0 & -1 & 2 
\end{pmatrix} \in \q M_{K+1}(\m R)
\]
Recall that $A$ is a positive definite matrix so that there exists $\lambda_K >0$ such that for all $\bm x \in \m R^{K+1}$, $\bm x^\top A \bm x \ge \lambda_K | \bm x|^2$. Therefore, 
\[ S \ge \lambda_K |q(\bm r)|^2 = \lambda_K |q_*(\bm y)|^2, \]
 and thus, with $\lambda$ being the minimum between $\lambda_k \frac{g_0\kappa }{2\alpha }$ and the implicit constant in the $O$ of \eqref{eq:est_dot_q*}, there holds
\[
\frac{\mathrm{d}}{\mathrm{d}t}\left[\frac{1}{q_*}\right]=-\frac{1}{q_*^{2}}\sum_{k=1}^{K-1}\dot{r}_{k}q'(r_k)\geq\lambda-\frac{1}{\lambda q_{*}} \left( \mathcal{N}^{2} + \sum_{k=1}^{K-1}r_k^{-1}q(r_k) \right). \]
This give the lower bound. For the upper bound, we simply notice that $S \le \| A \| |q_*(\bm y)|^2$, so that possibly decreasing $\lambda$ so that $\lambda < 1/\| A \|$, the conclusion follows as before.
\end{proof}

\bigskip
    
\emph{Bootstrap estimates.}
We introduce the following bootstrap estimates
\begin{equation}\label{BS-3}
    \mathcal N(t)  \le \delta^{3/4}  \qquad q_*(\bm y(t))  \le \delta^{3/2} \qquad \text{and} \quad
    b(t)  \le  \delta^{5/2}. 
	\end{equation}
	We set
	\begin{equation*}
	\TS=\sup\left\{ t\in [0,\infty)\hbox{ such that \eqref{BS-3} holds on $[0,t]$}\right\}\geq 0.
	\end{equation*}
We therefore emphasize that in the computations below, the implicit constants do not depend on  $\delta$.

\bigskip

\emph{Estimates on the damped components}

We argue as in the same step of the proof of Proposition \ref{pr:exist}, and we infer that for all $t \in [0,\TS]$,
\[ \q F(t), \q G(t) \lesssim \delta^2 (1+\delta t)^{-2} \]
so that
\begin{equation} \label{est:N_decay}
\q N(t) \lesssim \delta (1+\delta t)^{-1} \lesssim \delta,
\end{equation}
which strictly improves the bootstrap assumption on $\q N$.

\bigskip

\emph{Estimates on the distance}

We define the refined quantity
\[ D(t) := \frac{1}{q_*(\bm y(t))} e^{- 3 \lambda^{-1} \q M(t)} \]
Using \eqref{eq:M_lyap} and the lower bound in \eqref{le:q_* control}, we get
\[ \dot D(t) \ge \left(\lambda+\frac{2}{\lambda q_*}\left(\mathcal{N}^2-\sum_{k=1}^{K-1}r_k^{-1}q(r_k)\right)\right) \exp \left(-3 \lambda^{-1} \mathcal{M}\right).
\]

Now, for $t \in [0,\TS]$, $\exp \left(-3 \lambda^{-1} \mathcal{M}\right)=1+O\left(\mathcal{N}^2\right)=1+O\left(\delta^{2}\right)$, and $r_k \gtrsim |\ln \delta|$. Therefore, on $[0,\TS]$
\[
\dot D \ge \lambda (1 - O(\delta^{2})) - O(|\lambda \ln \delta|^{-1}) \ge \frac{\lambda}{2}.
\]
Integrating on $[0, t]$ we infer that $ D(t) \ge D(0)+\frac{\lambda}{2} t$. As we also have $D(0)^{-1} \le 2 q_*(\bm y(0)) \lesssim \delta^2$.
Therefore ,
\begin{equation} \label{est:q*_decay}
\frac{1}{q_*(\bm y(t))} \ge (1+ O(\delta^2) ) \left( \frac{2}{C\delta^{2}} + \frac{\lambda}{2} t \right) \ge  \frac{1}{C\delta^2} + \frac{\lambda}{4} t.
\end{equation}
In particular, $q_*(\bm y(t)) \lesssim \delta^2$, and this strictly improves the bootstrap estimate of $q_*(\bm y(t)$.

\bigskip

\emph{Transversality condition.}

Assume that there exists $t_* \in [0,\TS]$ such that $b(t) = \delta^{5/2}$. Then using \eqref{eq:dot_b}, 
\[ \dot b(t) \geq 2\nu^+ b(t) - C \delta^3 \geq 2\nu^+ \delta^{\frac{5}{2}} - C \delta^3 \geq \nu^+ \delta^{\frac{5}{2}} > 0. \]
Arguing as in the proof of Proposition \ref{pr:exist}, this transversality condition is enough to justify the existence of at least a point $\gh a^+_\sharp(0) \in \bar{\mathcal{B}}_{\mathbb{R}^K}(\delta^{\frac{5}{4}})$ such that $\TS = \infty$: we denote $\vec u_\sharp$ the associated solution via \eqref{def:ini2}, this is the sought for $K$ soliton solution. 

This closes the bootstrap.

\bigskip

\emph{Precise asymptotics of the solitons' centers of mass.}

The decays
\[ \q N, q_*(\bm y) = O(t^{-1}) \]
 follow from \eqref{est:N_decay} and \eqref{est:q*_decay}. As a consequence $r_k \ge \ln t + o(\ln t)$ for all $k=1, \dots, K-1$. It remains to compare $\q N$ and $q_*(\bm y)$.

For this, we consider the other refined quantity
\[ R(t) = \frac{1}{q_*(\bm y(t))} e^{\lambda^{-1} \q M(t)} \]
where $\lambda$ is given in Claim \ref{cl:dot_q*}: we use the upper bound in \eqref{le:q_* control}, we get that for large $t$, 
\[ \frac{d}{dt} \left[ \frac{1}{q_*(\bm y)} \right] \le \frac{1}{\lambda} + \frac{\q N^2}{\lambda q_*(\bm y)} + \frac{2}{\lambda \ln t} \le \frac{2}{\lambda} + \frac{\q N^2}{\lambda q_*(\bm y)}. \]
Therefore, using the above together with \eqref{eq:M_lyap} (and recalling that  $\q M = O(t^{-2})$)
\begin{align*}
\dot R(t) & = e^{\lambda^{-1} \q M(t)} \left( \frac{d}{dt} \left[ \frac{1}{\dot q_*(\bm y)} \right] + \frac{\dot M}{\lambda q_*(\bm y)} \right) \\
& \le e^{\lambda^{-1} \q M(t)} \left(  \frac{2}{\lambda} + \frac{\q N^2}{\lambda q_*(\bm y)} - \frac{\q N^2}{\lambda q_*(\bm y)}  + O(q_*(\bm y) ) \right) \le O(1).
\end{align*}
Integrating in time, this proves that  $R(t) \lesssim t$, and as a consequence (still using that $\q M$ is bounded), we infer that $q_*(\bm y(t)) \gtrsim t^{-1} \gtrsim \q N$. Therefore, the $O$ in \eqref{eq:ode_sys} can be replaced with solely $q_*(\bm y)^\theta$: system \eqref{eq:ode_sys} takes the form
\[ \dot y_k = \frac{1}{2\alpha} (g(y_{k} - y_{k-1}) - g(y_{k+1} - y_k)) + O(q_*(\bm y)^{\theta}). \]
Notice that at this point, the dynamics of the centers is essentially decoupled from the dynamics of the error $\vec \e$. 
Proposition \ref{prop:stab} now applies, and yields the desired asymptotics  \eqref{eq:zk} on $\bm y$ and so on $\bm z = \bm y + O(t^{-1})$.
\end{proof}

\begin{remark} \label{rem:last}
The analysis of this section is robust enough to prove a description result in the spirit of \cite{CMY} for global solutions with cylindrical symmetry, (actually it is really a stability result conditional to being near a multi-soliton configuration on a sequence of times). It reads as follows.

\medskip
Let $\vec u \in \mc C([0,+\infty), \q H_{\q C})$ be a multi-soliton in the sense of Definition \ref{def-multi wave}, with cylindrical symmetry. Then there exist $K \ge 1$ and $\sigma \in \{ \pm 1 \}$, such that $\vec u(t)$ admits the decomposition \eqref{eq:th:2}-\eqref{eq:th:z} for all large times.

\medskip

This result has two restrictions with respect to \cite[Theorem 1.1]{CMY}. One, it assumes cylindrical symmetry (which seems unavoidable at this point), and two, that all the nonlinear objects which appear in a decomposition are \emph{ground states}. In dimension 2 or higher, excited states could appear (even within cylindrical symmetry), so there is no ruling them out. We leave the proof to the interested reader.

\end{remark}

\appendix 

\section{Proof of Proposition \ref{prop:stab}} \label{app}

The proof is inspired by a similar analysis in \cite{MZ12}, see also \cite{CZ13}. 

For $k =1, \dots, K-1$, denote (as before) $r_k(t) = y_{k+1}(t) - y_k(t)$ . Let
\[  \underline r(t) = \min_{k=1, \dots, K-1} r_k(t)  \quad \text{and} \quad \bar r(t) = \max_{k=1, \dots, K-1} r_k(t), \]
and
\[ d(t) = \max_{k =1, \dots, K-1} r_k(t) -  \min_{k =1, \dots, K-1} r_k(t). \]
In this setting, we do not know yet that $r_k \ge \ln (t) + o(\ln t)$ (which holds in the context of the proof Proposition \ref{prop:1d constr}, as a consequence to the bootstrap argument); nevertheless, we still have a weaker  assumption which is $\underline r \to +\infty$.

Let us quantify the assumption: there exists $C, T>0$ such that for all $t \ge T$ and $k=1,\dots, K$, there holds
\[ \left| \dot y_k - \frac{1}{2\alpha} (g(r_{k-1}) - g(r_{k})) \right| \le C g(\underline r)^{\theta}. \]

\begin{lemma} \label{lem:d_bound}
There exists $C_1$ such that for all $t \ge 0$, $|d(t)| \le C_1$.
\end{lemma}

\begin{proof}
We argue by contradiction: if the conclusion does not hold, there exist two sequences of times $s_n < t_n$, $s_n \to +\infty$ and a sequence $A_n \ge n$ such that
\begin{equation} \label{def:s_n-t_n}
d(s_n) = A_n, \quad d(t_n) = A_n+3 \quad \text{and for all } t \in (s_n,t_n], \quad d(t) > A_n.
\end{equation}
Also, as $\underline r(s_n) \to +\infty$, we can assume that $\underline r(s_n) \ge 2$ and $g(\underline r(s_n)) \le 1$ for all $n$.

The argument rests on two claims \ref{cl:edo_1} and \ref{cl:edo_2} below which both rely on the study a suitable time rescaling of the ODE system. More precisely, consider
 the sequence of functions for $k=1,\dots, K-1$ and $n \in \m N$:
\[ h_{k,n}: s \mapsto r_k \left( s_n + \frac{s}{g(\underline r(s_n))}  \right) - \underline r(s_n). \]

Notice that for all $k=1,\dots, K-1$ and $s \ge 0$,
\begin{multline} \label{est:dot_hkn}
 \left| \dot h_{k,n} -  \frac{1}{2\alpha g(\underline r(s_n))} (- g(h_{k+1,n}+ \underline r(s_n)) + 2 g(h_{k,n}+ \underline r(s_n)) - g(h_{k-1,n}+ \underline r(s_n)) ) \right| \\
 \le 2C \frac{g(\underline r)^{\theta}}{g(\underline r(s_n))}.
 \end{multline}

\begin{claim} \label{cl:g}
There exists $M \ge 1$ such that for all $r \ge 2$ and $h \ge -1$, 
\[ 0 \le \frac{g(r+h)}{g(r)} \le M. \]
Furthermore, if $h \in \m R$ is fixed, then $\frac{g(r+h)}{g(r)} \to e^{-h}$ as $r \to +\infty$, with uniform convergence when $h$ remains bounded (or actually $o(r)$).
\end{claim}

\begin{proof}
Indeed, recall that due to \eqref{Qdec}, there exists $C_2 >0$ such that
\[  \forall y \ge 1, \quad \frac{1}{C_2 y^{\frac{d-1}{2}}} e^{-y}  \le g(y) \le \frac{C_2}{y^{\frac{d-1}{2}}} e^{-y}. \]
Then in the range specified for $r$ and $h$, 
\[ 0 \le \frac{g(r+h)}{g(r)} \le  \frac{C_2^2 }{\left( 1+ \frac{h}{r} \right)^{\frac{d-1}{2}}} e^{-h} \le  C_2^2 2^{\frac{d-1}{2}} e. \]
The limit $r \to +\infty$ for fixed $h$, is straightforward.
\end{proof}

\emph{Step 1. Limit system.}

Let $n \in \m N$. Define
\begin{gather*} s_n^* = \inf \{ s \ge 0 : \forall s' \in [0,s], \forall k=1, \dots, K-1, \ h_{k,n}(s') \ge -1  \},  \\ 
\text{and} \quad t_n^* = s_n + \frac{s_n^*}{g (\underline r(s_n))}. 
\end{gather*}
From Claim \eqref{cl:g}, we infer that \eqref{est:dot_hkn} now writes: for all $s \in [0,s_n^*)$,
\begin{multline} \label{est:dot_hkn2}
 \left| \dot h_{k,n} -  \frac{1}{2\alpha g(\underline r(s_n))} (- g(h_{k+1,n}+ \underline r(s_n)) + 2 g(h_{k,n}+ \underline r(s_n)) - g(h_{k-1,n}+ \underline r(s_n)) ) \right| \\
 \le 2CM g(\underline r(s_n))^{\theta-1}.
 \end{multline}
and so for all $s \in [0,s_n^*)$ (using again Claim \ref{cl:g} and $g(\underline r(s_n)) \le 1$),
\begin{equation} \label{est:equi_bd}
\forall s \in [0,s_n^*), \quad  |\dot h_{k,n}| \le \left(\frac{1}{\alpha} + 2C \right) M =: M_1.
\end{equation}

In particular, as $h_{k,n}(0) \ge 0$,  $s_n^* \ge 1/M_1$.

Fix $T \in (0,1/M_1)$ (independent of $n$) and let
\[ \hat t_n = s_n + \frac{T}{g(\underline r(s_n))} \le t_n. \]
In view of \eqref{est:equi_bd}, $|h_{k,n}(s)-h_{k,n}(0)|\le 1$ for all $s \in [0,T]$ and as by definition
\[
	h_{k,n}(s)-h_{k,n}(0)=r_{k}\left( s_n+\frac{s}{g(\underline r(s_n))}\right) -r_{k}(s_n),
	\]
we infer that 
\[ \forall t \in [s_n,\hat t_n], \quad |\underline r(t) -\underline r(s_n)| \le 1 \quad \text{and} \quad |\bar r(t) -\bar r(s_n)| \le 1. \]
As a consequence, $d(s) \le d(s_n)+2$ for all $s \in [s_n,\hat t_n]$ and so $\hat t_n < t_n$.

\bigskip

\emph{Step 2. Bound from above.}

\begin{claim} \label{cl:edo_1}
For all $\eta >0$, there exists $N = N(\eta)$ such that for all $n \ge N$,
\begin{equation} \label{est:bar_r}
\forall s \in [s_n,\hat t_n], \quad \bar r(s)  \le \bar r(s_n) + \eta.
\end{equation}
\end{claim}

\begin{proof}
Let $\bar h_n(s) = \max_{k} h_{k,n}(s)$. We claim that 
\begin{align} \label{est:barh}
\bar h_n(s) \le \bar h_n(0) + 2CM g(\underline r(s_n))^{\theta-1} s.
\end{align}

Indeed, if $s \in [0,T]$, let $j=j(s)$ be such that $\bar h_n(s) = h_{j,n}(s)$. As $g$ decreases on $[1,+\infty)$, we infer that
\[ \dot h_{j,n} \le 2CM g(\underline r(s_n))^{\theta-1}. \]
A continuity argument (and integration in time) allows to conclude to \eqref{est:barh}. As a consequence
\[ \forall s \in [0,T], \quad \bar h_n(s) \le \bar h(s) + 2CM T g(\underline r(s_n))^{\theta-1}. \]
Recall that $\underline r(s_n) \to +\infty$, so given $\eta >0$, define $N$ such that for all $n \ge N$,
\[ 2CM T g(\underline r(s_n))^{\theta-1} \le \eta. \]
As $\bar h_n(s) = \bar r\left( s_n + \ds\frac{s}{g(\underline r(s_n))} \right) - \underline r(s_n)$, the conclusion \eqref{est:bar_r} follows.
\end{proof}

\emph{Step 3. Bound from below}

\begin{claim} \label{cl:edo_2}
There exist $\delta_0 >0$ and $N_0 \in \m N$ such that for all $n \ge N_0$, there exists $\hat s_n \in [s_n, \hat t_n]$ such that
\[ \underline r(\hat s_n) \ge \underline r(s_n) + \delta_0. \]
\end{claim}

\begin{proof}
This is a consequence of Lyapunov functional on the limit system, which ensures that the corresponding $\underline r$ and $\bar r$ converge monotonically to one another.

Up to extracting a subsequence, we can assume that for all $k =1, \dots, K-1$, $$h_{k,n}(0) = r_k(s_n) - \underline r(s_n)$$ admits a limit $\ell_k \in [0,+\infty]$ as $n \to +\infty$.

Let us distinguish the indices $k$ depending on the values of this limit, and define
\begin{align*}
\q K_0 & = \left\{k \in \{ 1, \ldots, K-1 \} : \ell_k=0\right\}, \\
\q K_+ & = \left\{k \in \{ 1, \ldots, K-1 \} :  \ell_k \in (0,+\infty) \right\}, \\
\q K_\infty & = \left\{k \in \{ 1, \ldots, K-1 \} :  \ell_k=+\infty \right\}.
\end{align*}
Now, due to \eqref{est:equi_bd} and Ascoli's theorem, up to taking a subsequence, we can assume that for all $k=1, \dots, K-1$, $h_{k,n}$ converges uniformly on $[0,T]$ to either $+\infty$ if $k \in \q K_\infty$, or to a function $h_k \in \mc C([0,T])$ with $h_k(0) = \ell_k$ if $k \in \q K_0 \cup \q K_+$.

Let $\delta_1 = 1/2 \min(1, \ell_k: k \in \q K_+) >0$, and, using continuity of the $h_k$, let $S \in (0, T]$ be such that for all 
\begin{equation} \label{est:K+}
\forall k \in \q K_+ \cup \q K_\infty,  \forall s \in [0,S], \quad  h_k(s) \ge 2\delta_1.
\end{equation}

It remains to consider the case of the indices in $\q K_0$. Let $k \in \q K_0$, and consider the maximal interval of integers $\llbracket k_-,k_+ \rrbracket$ containing $k$ included in $\q K_0$, that is 
\[ \left\rrbracket k_-,k_+ \right\llbracket \subset \q K_0 \quad \text{and} \quad k_-,k_+ \in \q K_+ \cup \q K_\infty. \]
(We use our previous convention that $\ell_0 = \ell_K = +\infty$). Due to \eqref{est:dot_hkn}, taking the limit as $n \to +\infty$ (and recalling the limit in \eqref{cl:g}), we infer that for $j \in \left\rrbracket k_-,k_+ \right\llbracket$, $h_j \in \mc C^1([0,T])$ and
\[ \dot h_j =  \frac{1}{2\alpha} (- e^{-h_{j+1}} + 2 e^{-h_{j}} - e^{-h_{j-1}}) . \]
(and $h_j =0$, $h_{k_-}(0), h_{k+} (0) \in (0,+\infty]$).

This system was studied in \cite{MZ12,CZ13}, and it was shown that $h_j(s)>0$ for all $j \in \rrbracket k_-,k_+ \llbracket$ and $s > 0$ small, say $s \in (0, \tilde s_{k_-,k_+}]$ (see \cite[Claim 3.7]{MZ12}, and \cite[Proposition 2.1]{CZ13} for a Lyapunov functional), with $0 < \tilde s_{k_-,k_+} \le S$.

Let $\hat s = \min(S, \tilde s_{k_-,k_+} : k \in \q K_0)$ and 
\[ \delta_0 = \frac{1}{2} \min (\delta_1, \min (h_k(\hat s): k \in \q K_0)) >0, \]
then (recalling \eqref{est:K+}) we obtained that 
\[ \forall k \in \{ 1, \dots, K-1 \}, \quad h_k(\hat s) \ge 2 \delta_0. \]

By uniform convergence, there exists $N_0 \in \m N$ such that for all $n \ge N_0$ and for all $k \in  \{ 1, \dots, K-1 \}$, $h_{k,n}(\hat s) \ge \delta_0>0$. Letting $\ds \hat s_n := s_n + \frac{\hat s}{g(\underline r(s_n))}$, this translates to $r_k(\hat s_n) \ge \underline r(s_n) + \delta_0$, for all $k \in  \{ 1, \dots, K-1 \}$, or equivalently,
\[ \underline r(\hat s_n) \ge \underline r(s_n) + \delta_0. \qedhere \]
\end{proof}

\bigskip

\emph{Step 4. Conclusion}

Let $\delta_0, N_0$ be given by Claim \ref{cl:edo_2}. Then apply Claim \ref{cl:edo_1} with $\eta = \delta_0/2$, let $n = \max(N_0, N(\delta_0/2))$ (a time $s = \hat s_n$) so that
\[ \begin{cases}
\bar r(\hat s_n)  \le \bar r( s_n) + \delta_0/2 \quad \text{and} \\
\underline r(\hat s_n) \ge \underline r(s_n) + \delta_0.
\end{cases} \]
Taking the difference yields $d(\hat s_n) \le d(s_n) - \delta_0/2 = A_n - \delta_0/2$. But $\hat s_n \in [s_n,t_n]$, a contradiction to \eqref{def:s_n-t_n}.
\end{proof}

\begin{lemma} \label{lem:expansion_r}
For all $k=1, \dots, K-1$, $\ds r_k (t) = \ln t - \frac{d-1}{2} \ln \ln t + O(1)$.

(The implicit constant depends on $C_1, \alpha$ and $K$).
\end{lemma}

\begin{proof}

For $j =1, \dots, \lfloor K-1/2 \rfloor$, let 
\[ \Sigma_j = \sum_{k=j}^{K-j} r_k = r_j + \cdots +r_{K-j} \quad \text{and} \quad \Sigma = \sum_{j=1}^{K-1} 3^{1-j} \Sigma_j. \]
Recall that 
\[ 2 \alpha \dot r_k = -g(r_{k-1}) + 2 g(r_k) - g(r_{k+1}) + O(g(\underline r)^\theta). \]
Therefore, summing up, we get that for 
\[ \begin{cases}
2\alpha \dot \Sigma_j & = - g(r_{j-1}) + g(r_j) + g(r_{K-j}) - g(r_{K-j+1}) + O(g(\underline r)^\theta)  \\
& \hspace{60mm} \text{for } j =2,\dots, \lfloor K-1/2 \rfloor, \\
2 \alpha \dot \Sigma_1 & = g(r_1) + g(r_{K-1}) + O(g(\underline r)^\theta).
\end{cases} \]
As a consequence, we compute
\begin{equation} \label{eq:sigma}
2 \alpha  \dot \Sigma = \frac{2}{3} \sum_{j=1}^{\lfloor K-1/2 \rfloor} 3^{1-j} (g(r_{j}) + g(r_{K-j})) +  O(g(\underline r)^\theta).
\end{equation}
Now, due to Lemma \ref{lem:d_bound}, $|r_k - \underline r| \le C_1$ for all $k=1,\dots, K-1$. Therefore, summing up,
\begin{equation} \label{est:Sigma_r}
| \Sigma - \underline r|, |\Sigma - r_k| \le K C_1.
\end{equation}
These two bounds ensure that there exists $C_3 = C_3(K,C_1)$ such that
\begin{align*}
 \forall k=1, \dots, K-1, \quad \frac{1}{C_3} \le \frac{g(r_k)}{g(\underline r)} \le C_3, \quad \text{and} \quad 
 \frac{1}{C_3} \le \frac{g(\Sigma)}{g(\underline r)} \le C_3.
\end{align*}
Summing up, we conclude that there exists $C_4 = C_4(C_1,\alpha,K)$ such that
\[  \frac{1}{C_4} g(\Sigma)\le \dot \Sigma \le C_4 g(\Sigma). \]
Let $\Gamma$ be such that $\Gamma' = \frac{1}{g}$, notice that $\Gamma(r) \sim e^r r^{\frac{d-1}{2}}$ as $r \to +\infty$. Integrate in time of the previous bounds on $\dot \Sigma$: there exists $T \ge 1$ such that for all $t \ge T$,
\[ \frac{1}{2C_4} t \le \Gamma(\Sigma(t)) \le 2 C_4 t. \]
Arguing as for \eqref{equ:asyr} (specifically, we use the fact that if $\varphi (r) = e^r r^{\frac{d-1}{2}}$ and $\psi(r) = \ln r - \frac{d-1}{2} \ln \ln r$, both functions $\varphi$ and $\psi$ are increasing for $r \ge 1$, and for large $r$,    $r-1 \le \psi(\varphi(r)) \le r$), we conclude that
\[ \Sigma(t) = \ln t - \frac{d-1}{2} \ln \ln t + O(1). \]
Combining this with \eqref{est:Sigma_r} gives the conclusion.
\end{proof}

\begin{proof}[Conclusion of the proof of Proposition \ref{prop:stab}]
Recall that the assumption now writes
\[ \dot y_k = g(r_{k-1}) - g(r_k) + O(g(\underline r)^\theta). \]
From Lemma \ref{lem:expansion_r}, there hold that $\underline r(t) = \ln t - \frac{d-1}{2} \ln \ln t+ O(1)$ so that $g(\underline r(t))^\theta = O(t^{-\theta})$ is integrable in time. Now, if $\ds \underline y := \frac{1}{K} \sum_{k=1}^K y_k$, there hold
\[ \dot{\underline y} = O(g(\underline r)^\theta) = O(t^{-\theta}). \]
Hence $\underline y(t)$ admits a finite limit as $t \to +\infty$, and so is $O(1)$.

Now, using both Lemma \ref{lem:d_bound} and \ref{lem:expansion_r} again, there hold for any $k \in \{ 1, \dots, K \}$
\begin{equation} \label{est:yk-y1}
y_k - y_1 = r_{k-1} + \cdots + r_1 = (k-1) \ln t - \frac{d-1}{2} (k-1) \ln \ln t + O(1). 
 \end{equation}
Taking the mean in $k \in \{ 1 ,\dots, K \} $, we infer that
\[ \underline y - y_1 = \frac{K-1}{2} \ln t - \frac{d-1}{2} \frac{K-1}{2} \ln \ln t + O(1), \]
and so
\[ y_1 (t) =  - \frac{K-1}{2} \ln t + \frac{d-1}{2} \frac{K-1}{2} \ln \ln t + O(1). \]
Summing up with \eqref{est:yk-y1}, we obtain \eqref{expansion_yk}.
\end{proof}

\section*{Statements}

The authors confirm that they have no conflict of interest.

No data is associated with this research.

\end{document}